\documentclass[12pt]{amsart}
\usepackage{graphicx}
\usepackage{xcolor}
\usepackage{hyperref}
\usepackage{fullpage}
\usepackage{listings}

\def\zcomm#1{\hspace*{2pc}\code{#1}}
\def\zoik#1#2{\hspace*{4pc}\hyperlink{#1}{#2}\\}
\def\zoikt#1#2{\hyperlink{#1}{#2}}

\def\xoik#1#2{\hspace*{0pc}\hyperlink{#1}{#2}\\}

\def\4{4D}
\def\3{3D}
\def\2{2D}
\def\1{1D}

\def\kp{KnotPlot}
\def\pui#1{\texttt{#1}}   
\def\panel#1{\texttt{#1}}
\def\code#1{\texttt{#1}}
\def\comm#1{\hyperlink{#1}{\commandname{#1}}}   
\def\demo#1{\texttt{#1}}
\def\param#1{\hyperlink{#1}{\commandname{#1}}} 
\def\eric#1{}

\newcommand\commandname[1]{\texttt{#1}}
\newcommand\command[1]{\hspace*{0.2in}\texttt{#1}}
\newcommand\tabname[1]{\texttt{#1}}
\newcommand\roller[1]{\texttt{#1}}
\newcommand\button[1]{\texttt{#1}}

\newdimen\ppplen
\def\ppp#1#2{\vtop{\hbox to \ppplen{\includegraphics[height=\ppplen]{#1}}\hbox to \ppplen{\hss#2\hss}}}
\newdimen\qqqwidth
\newdimen\qqqheight
\def\qqq#1#2{\vtop{\hbox to \qqqwidth{\includegraphics[height=\qqqheight]{#1}}\hbox to \qqqwidth{\hss#2\hss}}}
\def\qqq#1#2{\vtop{\hbox to \qqqwidth{\includegraphics[width=\qqqwidth]{#1}}\hbox to \qqqwidth{\hss#2\hss}}}
\begin{document}

\title{An Introduction to KnotPlot}


\author{Robert G. Scharein}
\address{Hypnagogic Software, Vancouver, BC, V6K 1V6, Canada}
\email{rob@hypnagogic.net}
\urladdr{www.knotplot.com}

\author{Eric J. Rawdon}
\address{Department of Mathematics, University of St.~Thomas, Saint Paul, MN 55105, USA}
\email{ejrawdon@stthomas.edu}
\urladdr{george.math.stthomas.edu/$\sim$rawdon}


\begin{abstract}
We give a brief introduction to the software \kp.
The goals of this chapter are twofold: 1) to help a new user get
started with using \kp{} and 2) to provide
veteran users with additional background and functionality
available in the software.
\end{abstract}

\def\zoik#1#2{\hspace*{4pc}\hyperlink{#1}{#2}\\}

\maketitle

\section{Introduction}

\kp{} is a program for visualizing and interacting with 3D
mathematical knots.  It was created in 1992 by RGS
in the Department of Computer Science at the University of British
Columbia \cite{Scharein_1998}.  Since then the program has steadily
accreted new features and has been applied to problems in mathematical
knot theory such as the minimal edge number of knots for various
cases \cite{Scharein_1998,RawSch,Scharein_2009,Ishihara_2012} and the
visualization of topology changes in
DNA \cite{TopoICE-X,TopoICE-R,TopologicalComplexity}.


This document is not intended as a comprehensive guide to \kp, but as
an introduction to aid new users to start using the software to
perform useful experiments.  Please refer to the manual for more
information:  A hyperlinked HTML manual is available at \cite{KPMhtml} and a
PDF version can be downloaded from \cite{KPM}.  The version of
KnotPlot described here is the current one at the time of writing
(March 2023).  For consistency with what is in this document, make
sure your copy of \kp{} was compiled no earlier than this date.  You
can find the compile date using the command \comm{version}.

\eric{The not linear assumes the reader is reading the PDF.  Not sure how likely this is.}

This chapter is not linear.  Each section is generally self contained.  The list below provides a short description of what is in each section.
We highlight three particular sections:
\zoikt{commands}{Section \ref{sec:commands}} that lists and describes commonly used commands,
\zoikt{parameterssection}{Section \ref{sec:params}} that shows the user how to change parameter
values (like the background color), and
\zoikt{ActivitySection}{Section \ref{sec:list}} that lists commonly used commands in popular
user activities.

\eric{ I agree it looks awful, how about some space between each line?   } 

{\bfseries

\xoik{setupsection}{\ref{sec:setup} -- Setting it up:\\\hspace*{0.19in} Downloading and setting up \kp\vspace*{0.1in}}
\xoik{LoadNSave}{\ref{sec:loadsave} -- Loading and saving:\\\hspace*{0.19in} Basic input/output for knot/link models\vspace*{0.1in}}  
\xoik{changing}{\ref{sec:view} -- Changing the view or embedding:\\\hspace*{0.19in} Rotating a model\vspace*{0.1in}}
\xoik{simple}{\ref{sec:simple} -- Relaxing knots and links examples:\\\hspace*{0.19in} How to relax a knot/link so that it is visually appealing\vspace*{0.1in}}
\xoik{pictures}{\ref{sec:pict} -- Making pictures:\\\hspace*{0.19in} Directions for making images of 3D knots/links\vspace*{0.1in}}
\xoik{DiagramsSection}{\ref{sec:eps} -- Making knot diagrams:\\\hspace*{0.19in} Directions for making different types of knot diagram images\vspace*{0.1in}}
\xoik{creating}{\ref{sec:creating} --  Creating knots and links:\\\hspace*{0.19in} Sketching, editing, and constructions of special classes of knots/links\vspace*{0.1in}}
\xoik{commands}{\ref{sec:commands} --  Useful commands:\\\hspace*{0.19in} List and description of highly-used \kp{} commands\vspace*{0.1in}}
\xoik{parameterssection}{\ref{sec:params} -- Parameter values:\\\hspace*{0.19in} Guide to some parameter values and how to change them\vspace*{0.1in}}
\xoik{uisection}{\ref{sec:ui} -- User interface:\\\hspace*{0.19in} Guide for windows, tabs, boxes, and rollers seen when launching \kp\vspace*{0.1in}}
\xoik{aloadsave}{\ref{sec:aloadsave} -- KnotPlot distribution:\\\hspace*{0.19in} Information about where models are loaded from and save to, as well as a guide to directories in the \kp{} installation\vspace*{0.1in}}
\xoik{advanceddynamics}{\ref{sec:dynamics} -- Advanced dynamics:\\\hspace*{0.19in} How to change the dynamics for relaxing knots/links\vspace*{0.1in}}
\xoik{scripting}{\ref{sec:scripting} -- Scripting and running without graphics:\\\hspace*{0.19in} Using files of commands to perform experiments using \kp\vspace*{0.1in}}
\xoik{TabsSection}{\ref{sec:tabs} -- Tabs:\\\hspace*{0.19in} Guide to the tabs which we do not cover in detail elsewhere\vspace*{0.1in}}
\xoik{ActivitySection}{\ref{sec:list} --  Commands listed by activity:\\\hspace*{0.19in} Commands grouped by common user activities}
}

\hypertarget{setupsection}{}

\section{Setting it up}
\label{sec:setup}

KnotPlot may be downloaded from \cite{KPdownload}.  It is available
for macOS (both Intel and ARM architectures), Windows 10/11, and
multiple flavours of Linux.  Instructions are also given on the
download page for setting up project directories.  This is important
for conducting research with KnotPlot in order to keep individual
projects self contained.

\begin{figure}[htb] 
\hbox to \textwidth {\hss\includegraphics[width=\textwidth]{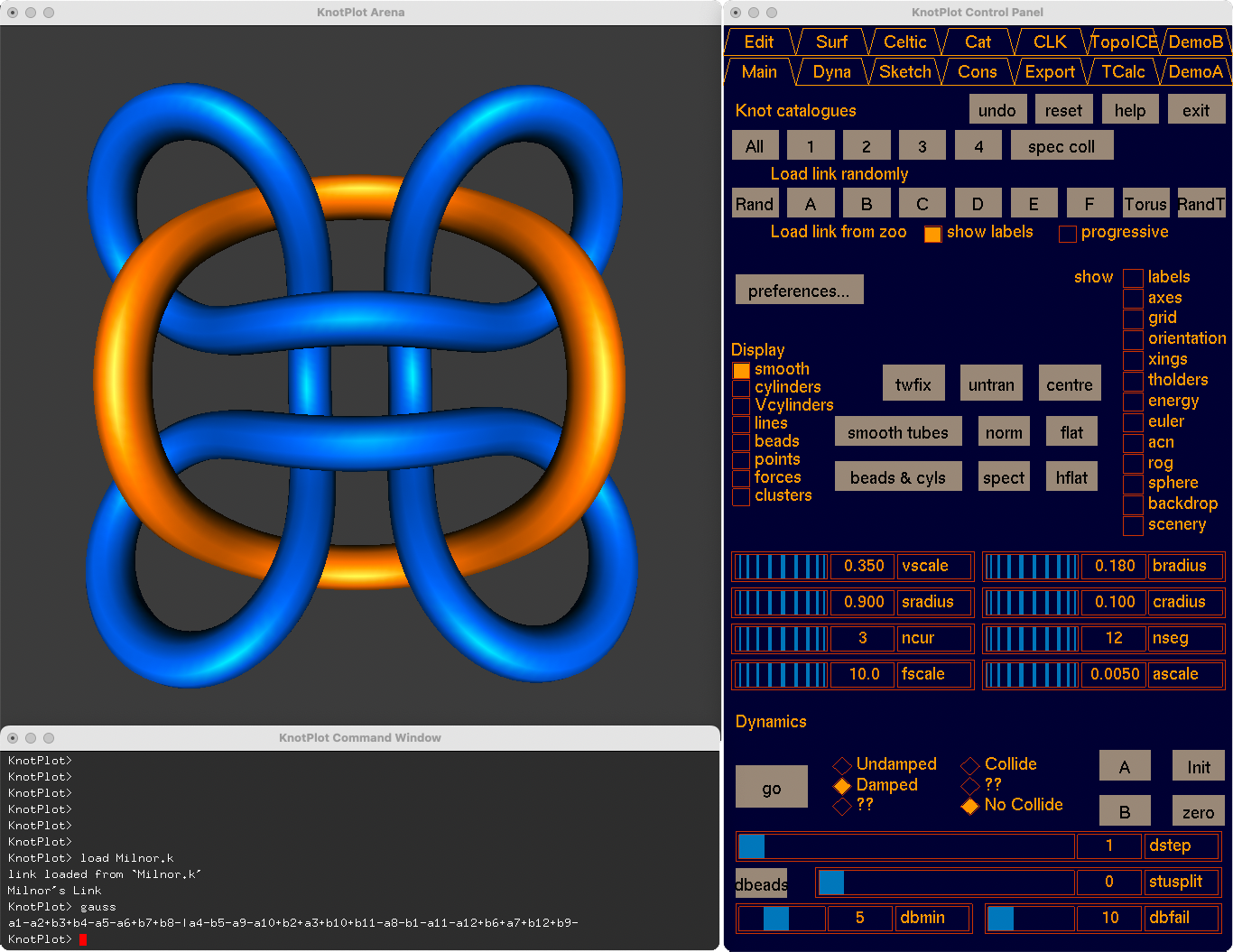}\hss}
\caption{\kp{} showing the Arena (upper left), Control Panel (right, set to the \panel{Main} tab) and Command Window (bottom).}
\label{fig:KnotPlot}
\end{figure}
 
Figure \ref{fig:KnotPlot} shows what KnotPlot looks like upon startup,
after a knot of interest is loaded.  KnotPlot opens three windows:
\emph{KnotPlot Arena}, \emph{KnotPlot Control Panel}, and
\emph{KnotPlot Command Window}.  The knot is displayed in the
\emph{KnotPlot Arena} where interactions using mouse and keyboard can
be used.  The \emph{KnotPlot Control Panel} is used for many tasks.
It is divided into a number of sub-panels that may be accessed using
one of the 14 tabs at the top of the panel.  In this document when we
refer to a specific sub-panel we use the names shown, for example
\panel{Main}, \panel{Sketch}, and \panel{DemoA}.  We discuss a few of
these tabs in detail, please refer to Section
\hyperlink{TabsSection}{\ref{sec:tabs}} to see a brief description of
the rest.  The final option for interacting with KnotPlot is to use
the \emph{KnotPlot Command Window}.  Here we see that that Milnor's
Link has been loaded and its extended Gauss code calculated by
entering the commands \code{load Milnor.k} and \comm{gauss}.

\eric{changed width of image to textwidth.  those two knots are in the distribution that people will download}

At the moment, the user interface looks identical on the various
platforms, having been implemented in OpenGL \cite{OpenGL} and GLUT
\cite{GLUT}.  This ``Classic KnotPlot'' does not use the native buttons,
sliders, and menus of the platforms they are running on.  For example,
KnotPlot's menus are accessed by right clicking on the Control Panel.
While this universal look has a certain appeal, it may make KnotPlot
more difficult to learn than it should be.  Currently in beta testing
is KnotPlot Redux, a cross-platform re-implementation of the user
interface with native look and feel, with ``Classic'' mode available
as a legacy option.

\hypertarget{LoadNSave}{}

\section{Loading and saving}
\label{sec:loadsave}

We recommend that users create a project directory (folder)
and start \kp{} with that directory set as the working directory.
If the working directory contains user's knot files, they may be loaded with something like\\
\zcomm{load myknot}\\
where \texttt{myknot} is the name of knot file.
If the file specified does not exist in the working directory, \kp{} will search a list of directories to find the file.
For more information on this, please refer to the \comm{path} command. 

KnotPlot comes with the entire Rolfsen Appendix C catalogue of knots and links \cite{Rolfsen}.
The KnotPlot names are formed in a simple way, for example\\
\zcomm{load 3.1}\\
loads the trefoil $3_1$ and\\
\zcomm{load 9.2.14}\\
the link $9^2_{14}$.

Most of the thousands of knots and links in \kp{}'s various catalogues are stored in 
a compact binary format, but users can use a simple plain text format to load a configuration.
This is one vertex per line with coordinates separated by spaces, and a blank line separating different components.
For  example the following represents the Hopf link $2^2_{1}$.\\
\zcomm{0.10 -3.29 -0.49}\\
\zcomm{1.11 0.69 0.13}\\
\zcomm{-1.64 -0.27 0.26}\\
\zcomm{}\\
\zcomm{0.01 2.30 0.44}\\
\zcomm{-0.07 2.10 -0.87}\\
\zcomm{0.48 -1.55 0.53}\\

Saving knots is done with the \comm{save} command,\\
\zcomm{save myknot}\\
that will save to the working directory.
KnotPlot supports saving in many different file formats.
By default KnotPlot saves in its compact binary format.
Appending a file extension to the file name will save in a different format, as in\\
\zcomm{save something.txt}\\
for simple plain text and\\
\zcomm{save something.vect} \\
 for Geomview's VECT format \cite{VECT}.

\hypertarget{changing}{}

\section{Changing the view or embedding}

\label{sec:view}

The Arena shows the currently loaded knot using a \emph{virtual
camera} located at the starting position of $(0, 0, 8)$ with the centre of view
being the origin.  A viewing transformation is applied to the knot and
may be modified using the following mouse or keyboard controls:
\begin{itemize}
\item Left click and drag rotates the view.  You can also constrain
the rotation using the keyboard.  Press and hold one of \code{x},
\code{y}, or \code{z} before left clicking causes the rotation to be
about one of those coordinate axes.  These axes are in the coordinate
space of the knot, which may be already rotated.  If you want to
rotate about axes fixed to the screen, press and hold  \code{i}, \code{j}, or
\code{k} instead (\code{k} goes into the screen).  For precise control
over rotation, use the \comm{rotate} command.

\item Middle click and drag to translate the view (or press the
option/alt key (or ctrl key on some systems) and left click, and then drag).

\item Right click and drag to scale the view (or press the shift key
and left click, and then drag).  The viewing scale can be set to an
exact value using the \param{vscale} parameter.
\end{itemize}

These actions change only the view and not the actual embedding of the knot. 
Another command, \comm{about}, does the complement, namely changes the
actual embedding but not the view. 
The two may look exactly the same, for example compare \code{rotate x 33} to \code{about x 33},
but they are not equivalent. 
The \demo{rotate} demo in the \texttt{Tutorials} section on  \panel{DemoA} 
gives a good introduction to the differences
between \comm{rotate} and \comm{about}.  
The \hyperlink{rotatefix}{\commandname{rotate fix}} command is useful for transferring a view rotation to a rotation of the embedding.

\hypertarget{simple}{}

\section{Relaxing knots and links examples}

\label{sec:simple}

Sometimes we might have a knot in a rather messy conformation and we would just like to get a nice picture of it.
We can use the controls on \panel{Main} to accomplish this.
Try the following: 
\begin{itemize}

\item If necessary click on the \pui{reset} button to get \kp{} into its starting configuration. 

\item Load a knot of interest, say the trefoil\\
\zcomm{load 3.1}\\
or load it from the Knot Zoo using the \pui{A} button on \panel{Main}.

\item Click on the \pui{beads \& cyls} button so that you can see edges and vertices of the knot. 

\item The \pui{dbeads} button deletes edges in the knot without changing the knot type. 
Click on this button several times until nothing more can be deleted.

\item Click the \comm{go} button to start relaxing.   That button will turn green to indicate relaxation is happening. 
You should see that the knot gets stuck. 

\item Now increase the value of \param{stusplit} using the slider.   This allows any stuck edges to be split, in a way that  maintains knot type. 
The knot should now look much more trefoilesque. 

\item Keep the relaxation going and get an even nicer trefoil by entering the \comm{split} into the Command Window. 
This splits all the edges.  You should see something looking much like what you started with. 

\item If you feel like trying another example, just keep the relaxation going and load a knot or link.

\end{itemize}

Here is a real-world example which also handles rotating the knot so that it can be loaded
later from the given position.
We will use a knot that comes from an experiment that RGS did to generate ``compact'' versions of knots fitting into a spherical cavity.
Load one of these with\\
\zcomm{load compact.k}\\
and apply the method above, relaxing with a non-zero value of \param{stusplit} and repeatedly clicking on \pui{dbeads}.
When you have simplified it to your liking, stop the relaxation and click on the \pui{centre} button.
It will likely not be in an orientation showing a minimal crossing number.
Try rotating it until it looks as minimal as possible. 
Now what we would like to do is to save the \emph{rotated} version.  We need to apply the viewing transform to change the coordinates of the knot.
For this we use the \hyperlink{rotatefix}{\code{rotate fix}} command
and then \hyperlink{save}{\code{save transformedknot}}.  At this point,
if you exit \kp{}, restart \kp{}, and then \hyperlink{load}{\code{load transformedknot}}, you will see the view of the knot as you saved it (unless you
scaled the knot using \hyperlink{vscale}{\code{vscale}}).

\begin{figure}[htb] 
\hbox to \textwidth {\hss\includegraphics[height=5cm]{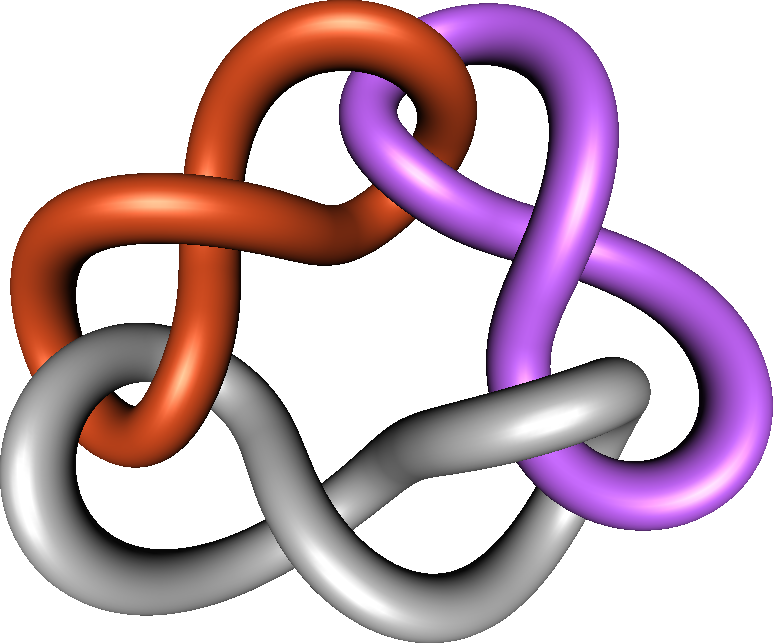}\hss}
\caption{The link $9^3_7$.}
\label{fig:3complink}
\end{figure}

Here is another example.  This creates a link with many vertices and then
simplifies the link.\\
\zcomm{load 9.3.7}\\
\zcomm{(press the beads \& cyls button)}\\
\zcomm{\#2 split \% splits each edge into two pieces, and does this twice}\\
\zcomm{fitto mindist 0.5 \% scales so min dist between beads is 0.5}\\
\zcomm{go beadlimit 4000 \% increase this value if there are many edges}\\
\zcomm{(press the go button, it will turn green)}\\
\zcomm{(press the dbeads button repeatedly)}\\
\zcomm{stusplit = 3}\\
After a few seconds, the link should look similar to Figure \ref{fig:3complink}.
Press the \texttt{go} button again to stop the relaxation.

\begin{figure}
	\ppplen 0.333\textwidth \hbox to \textwidth{
	       	\ppp{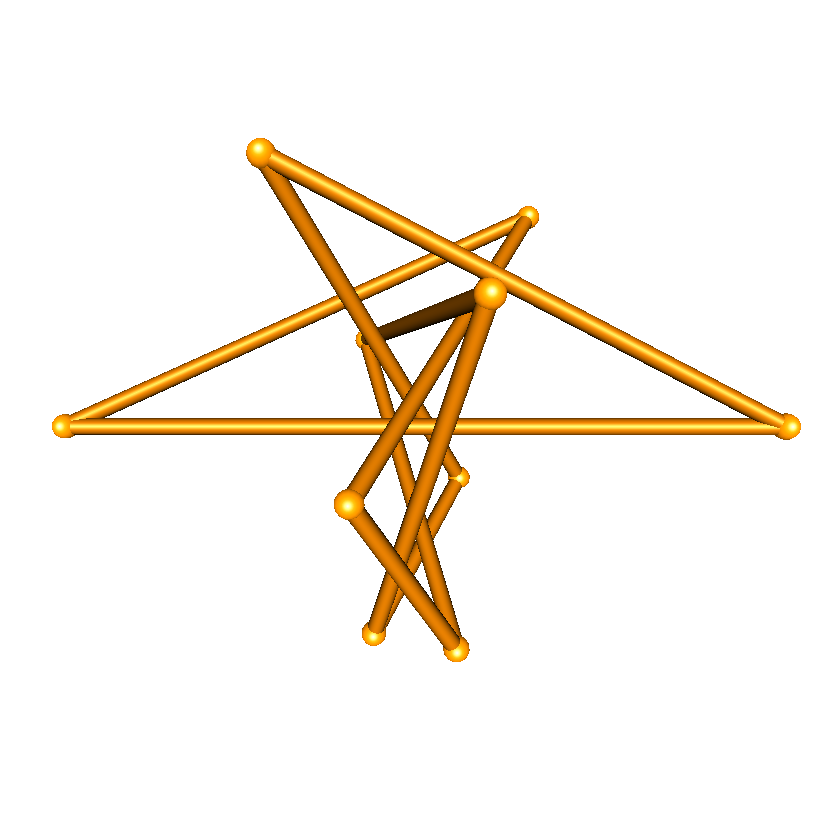}{\code{mode cb}}\hss
		\ppp{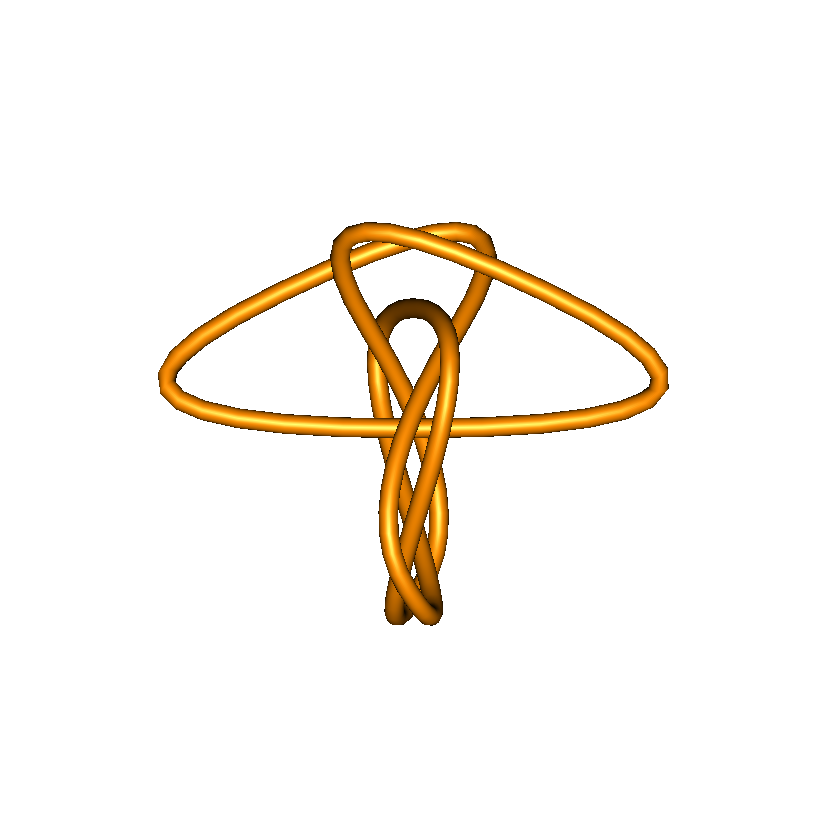}{\code{mode s}}\hss
		\ppp{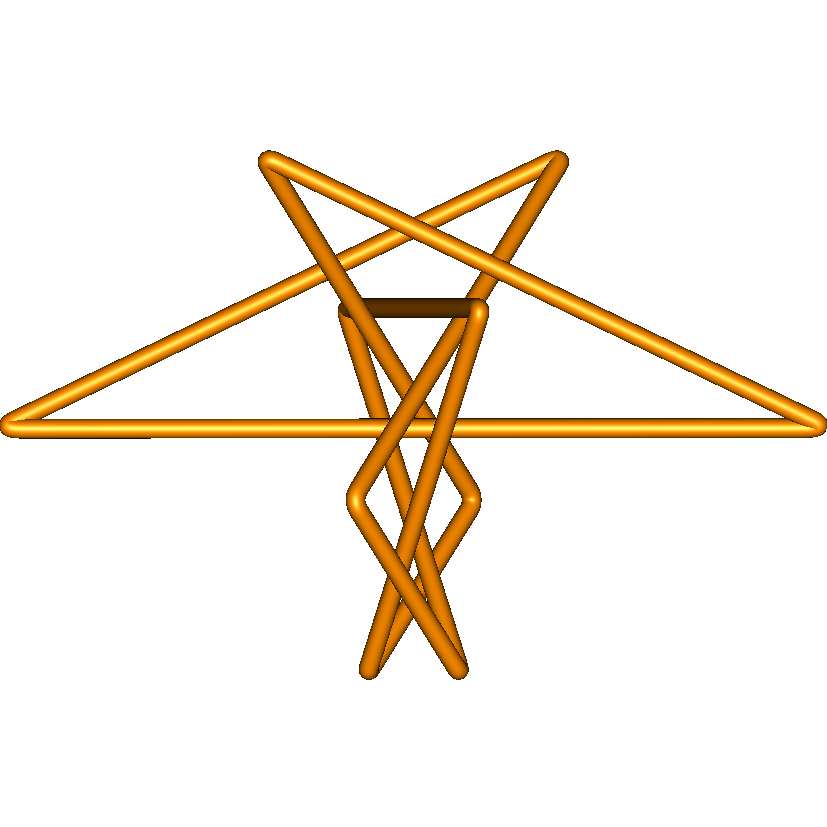}{\code{mode s} after splitting}}
  \caption{Various renderings of minimal stick $10_{124}$.}
  \label{fig:PLalt}
\end{figure}

\hypertarget{pictures}{}

\section{Making pictures}
\label{sec:pict}
Making full colour images of what is seen in the KnotPlot Arena can be made with the \comm{imgout} command.  
The colour of different components may be changed in various ways using the \comm{colour} command. 
Also, for publication purposes it is often desirable to have a white background rather than KnotPlot's grey colour. 
Entering the sequence of commands\\
\zcomm{load 9.3.7}\\
\zcomm{luxo}\\
\zcomm{background = white}\\
\zcomm{colour 0 grey}\\
\zcomm{colour 1 rgb:0.85/0.3/0.12}\\
\zcomm{colour 2 rgbi:199/102/255}\\
\zcomm{imgout 3complink}\\
results in a file called \texttt{3complink.png}, which is shown in Figure \ref{fig:3complink}.
The \comm{luxo} command is just a shortcut to setting the values of \param{ncur} and \param{nseg}. 

The next example is the kind of thing that we sometimes need to do to
get a good illustration.  Consider the knot on the left in Figure
\ref{fig:PLalt}.  This is a minimal stick representative of a knot
found in RGS's survey \cite{Scharein_1998}, with an embedding that
results from minimizing Simon's minimum distance energy
\cite{SimonMD}.\footnote{Note that the terms ``minimal'' and
``minimize'' are used in three different senses here.} It was then
aligned using \hyperlink{alignaxes}{\texttt{align axes}}.  The knot
exhibits a beautiful symmetry along all three axes but this is not
completely evident in the left-hand figure because we are viewing in
perspective.  Suppose we want to show the symmetry and also draw the
knot as just a polygon, without the little spheres showing the
vertices.  Switching to orthographic view with the command
\code{orthographic} and to smooth mode with command \code{mode s} we get the
middle figure, which may not be ideal.  This is because the
smooth tubes are drawn by using B\'ezier splines that interpolate the
edges.  The interpolation gives a correct look for knots only if the
turning angle is not too large between edges.  To get around this problem, we
use the trick of splitting edges multiple times to get the view on the
right.  The commands that generated Figure \ref{fig:PLalt} are below.  \\
\zcomm{luxo;background = white; mode cb; vscale = 0.8 } \\
\zcomm{sradius = 0.12; cradius = 0.12}\\
\zcomm{load ms/10.124; about z -90}\\
\zcomm{display true;imgout polygonalLeft }\\
\zcomm{mode s;ncur = 11; ortho}\\
\zcomm{display true; imgout polygonalMiddle }\\
\zcomm{\#4 split \ \ \ \ \% repeat split 4 times}\\
\zcomm{display true; imgout polygonalRight }\\

Note that several commands can be typed on the same line using
semicolons for separators
and that \code{\%} is the KnotPlot comment character (just like \LaTeX).
The command
\hyperlink{dispt}{\commandname{display true}} command is often used in
scripts run from the command window to force a redraw of the frame
buffer at that point.  Otherwise the frame buffer does not get redrawn
until the script finishes.  If you type these commands in manually,
you do not need \code{display true}.  Note that these commands can be
saved to a file, say \texttt{myscript.kps}, and then loaded through
the command window via \code{< myscript.kps}.  There is also a no
graphics mode to \kp{} which we discuss in \zoikt{scripting}{Section \ref{sec:scripting}}.

\begin{figure}
	\ppplen 0.333\textwidth   \hbox to \textwidth{\hss
	       	\ppp{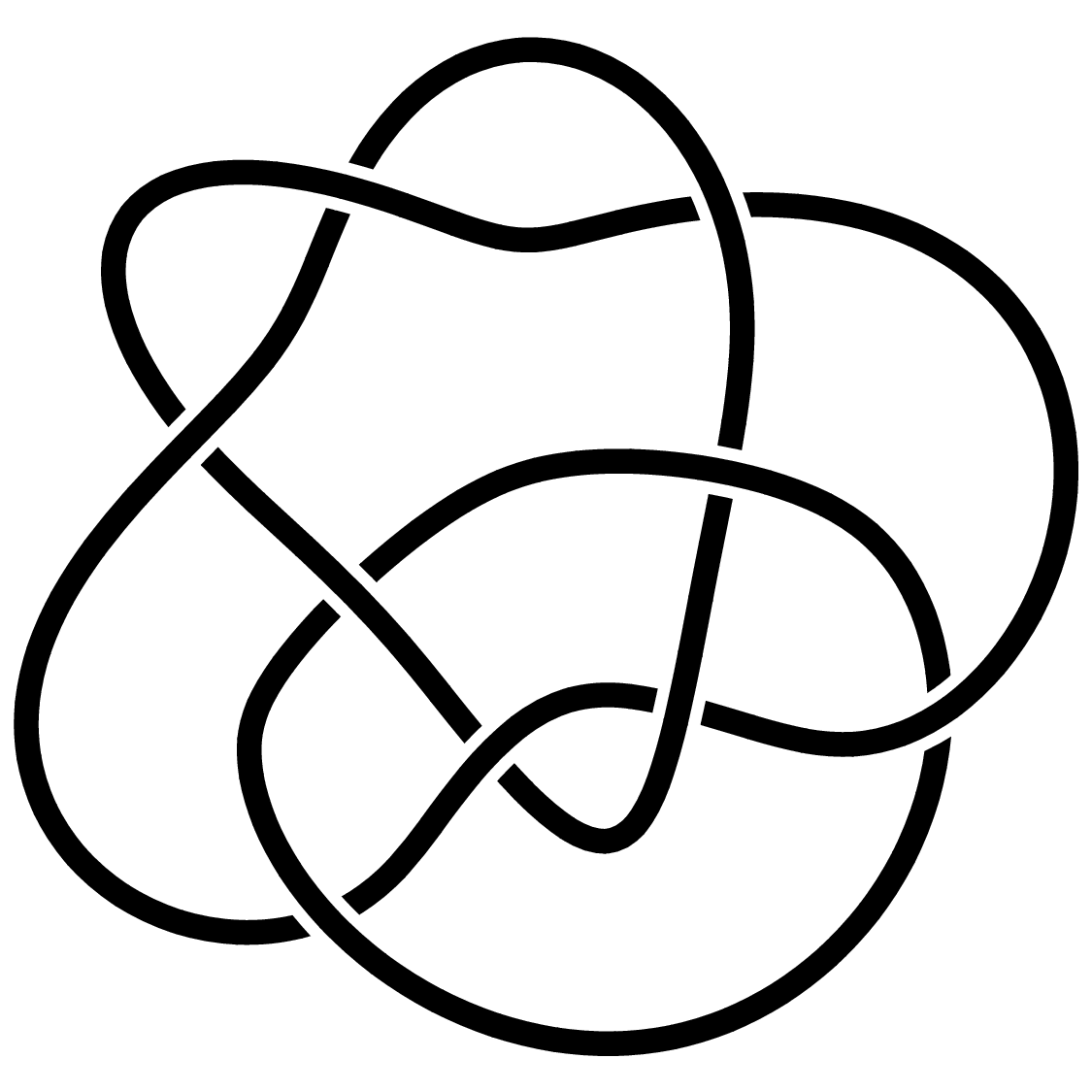}{\code{psmode = 40}}\hss
		\ppp{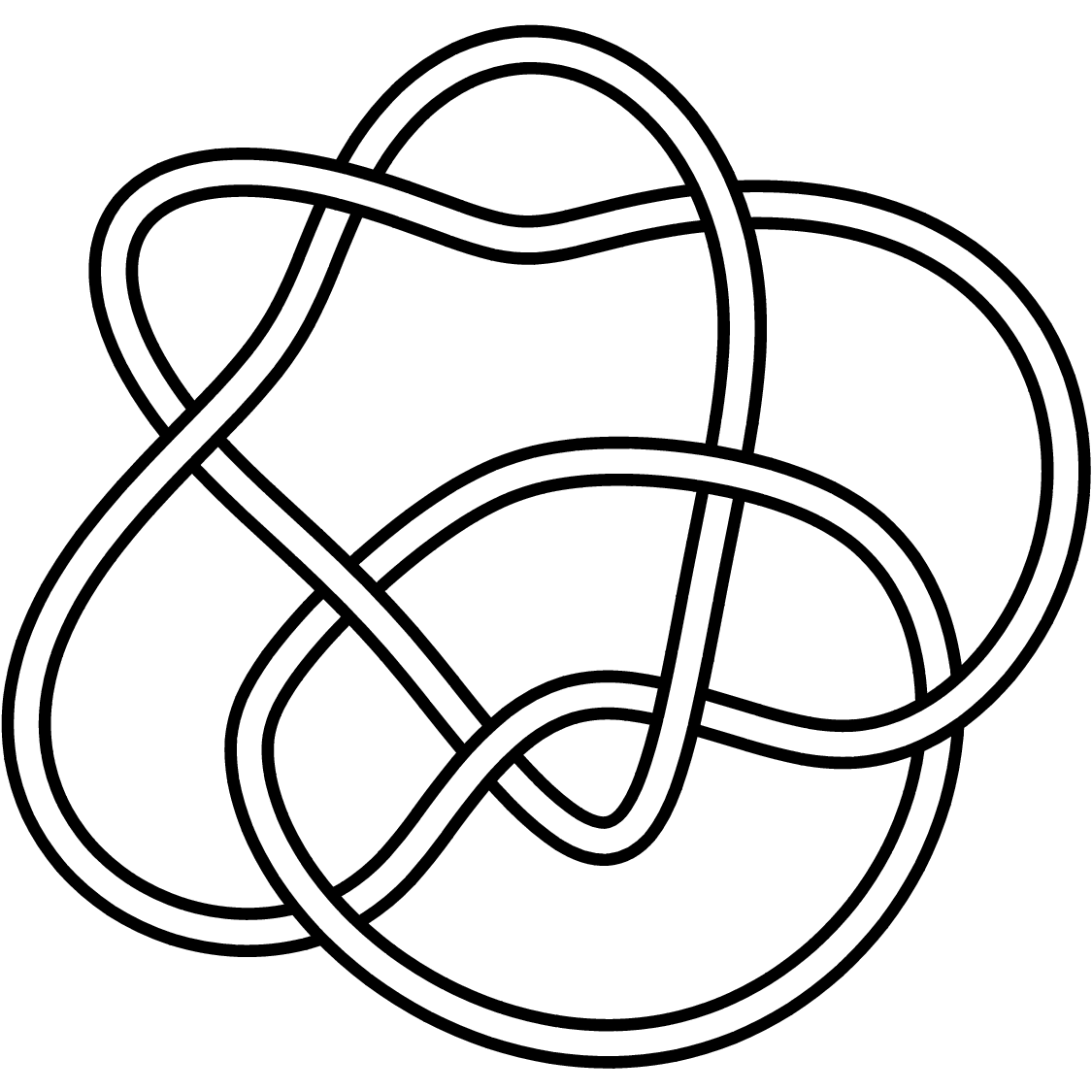}{\code{psmode = 41}}\hss
		\ppp{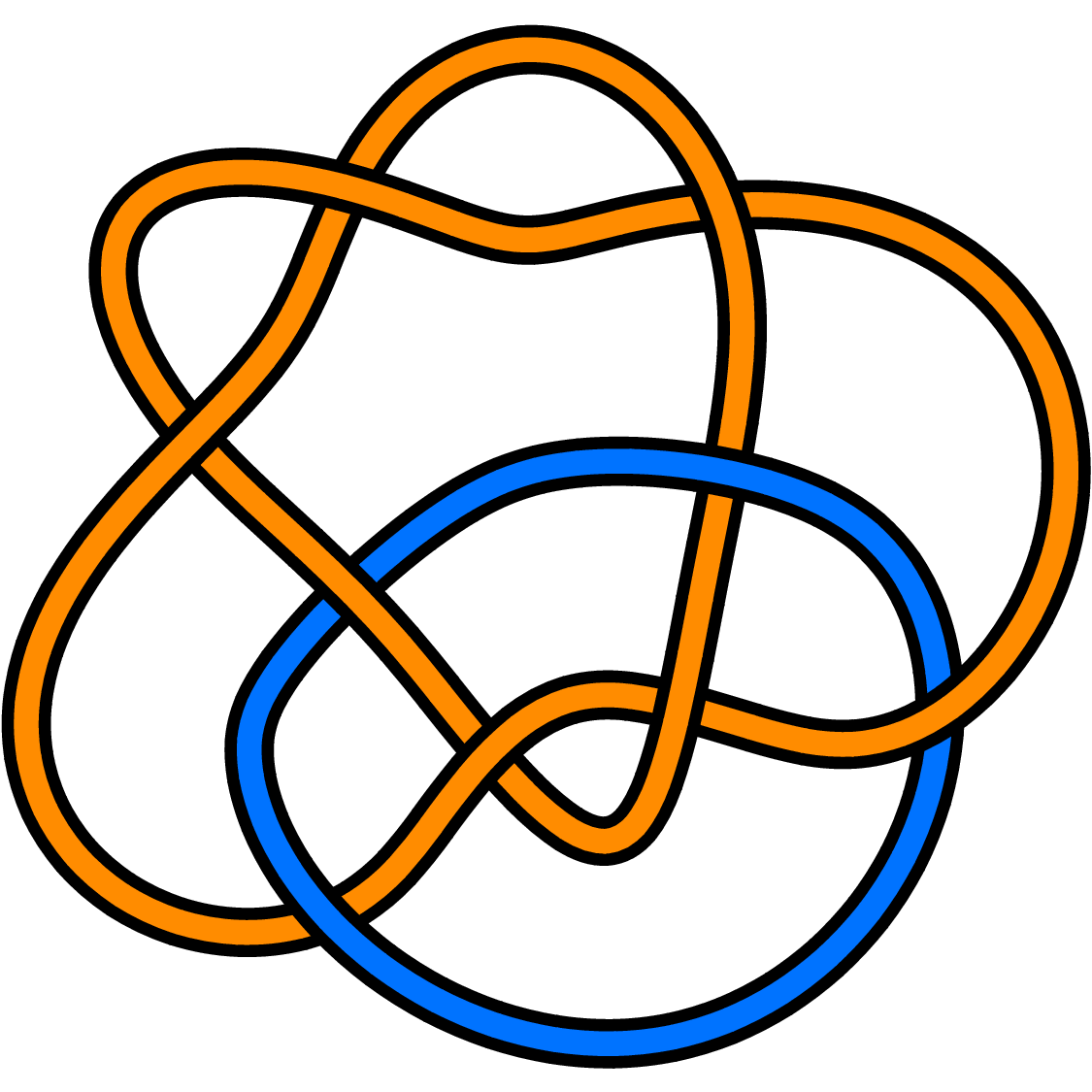}{\code{psmode = 45}}\hss}
  \caption{The link $9^2_{37}$ in different styles.}
  \label{fig:EPS}
\end{figure}

\begin{figure}
	\ppplen 0.333\textwidth    \hbox to \textwidth{\hss
	       	\ppp{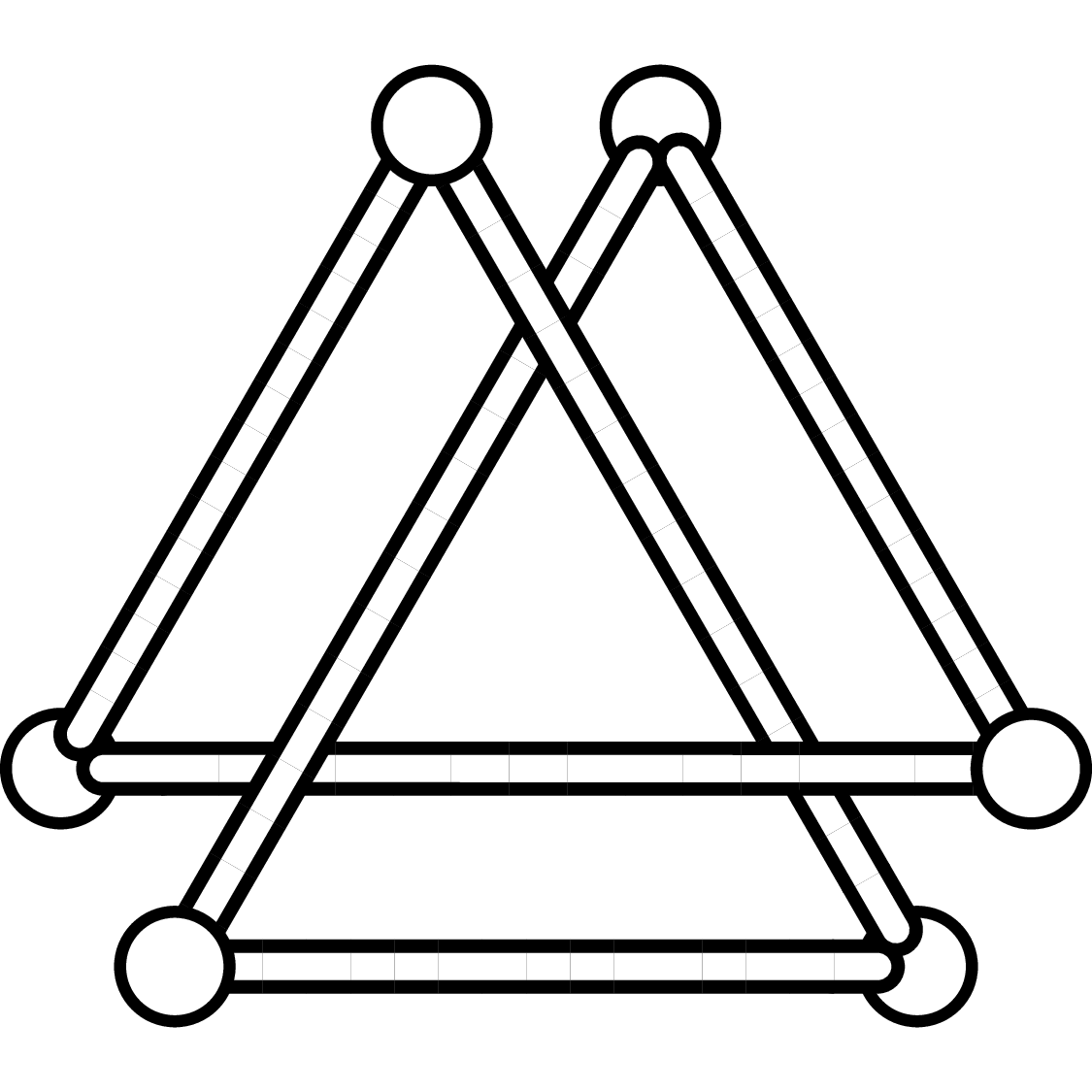}{Minimal stick}\hss
		\ppp{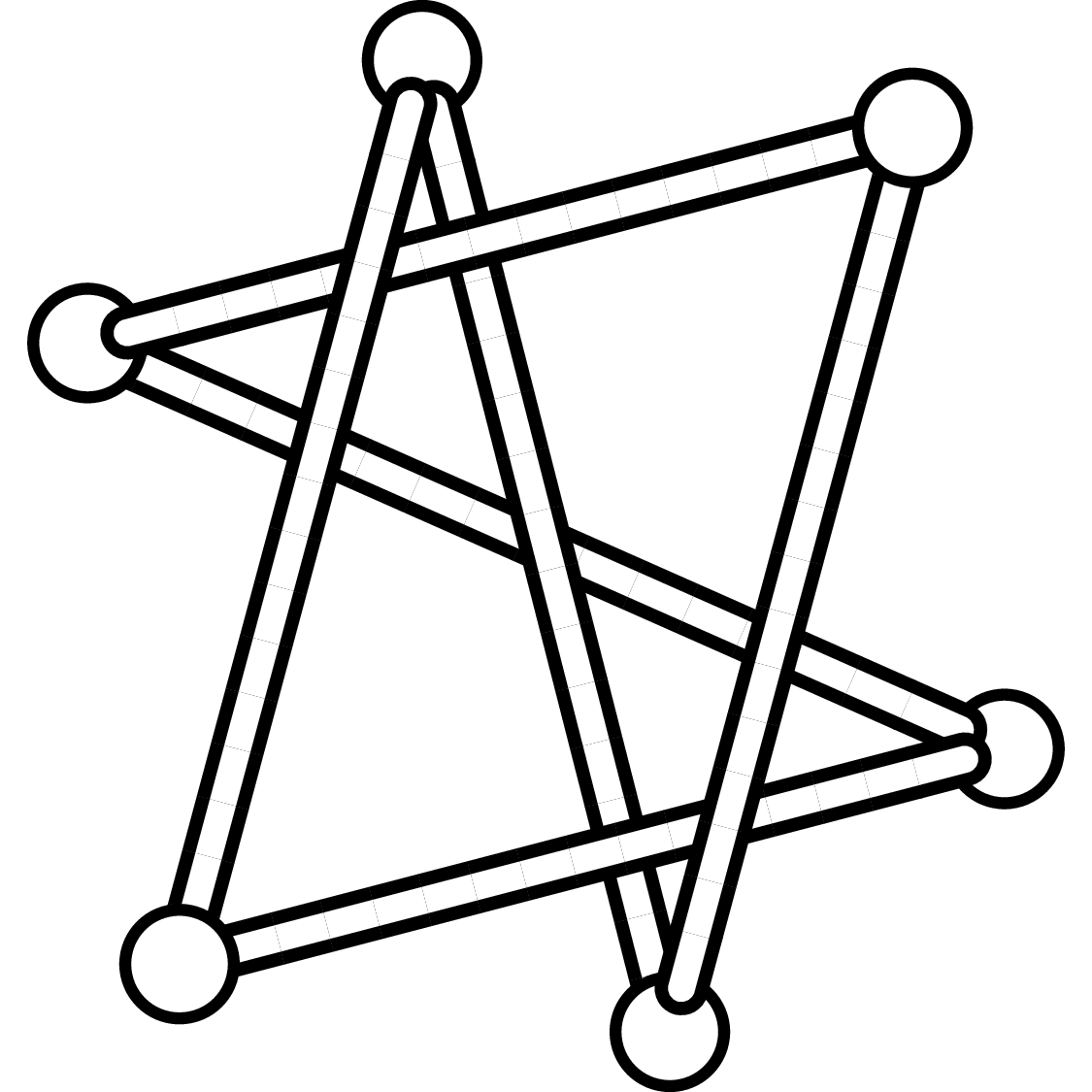}{Equilateral minimal stick}\hss
		\ppp{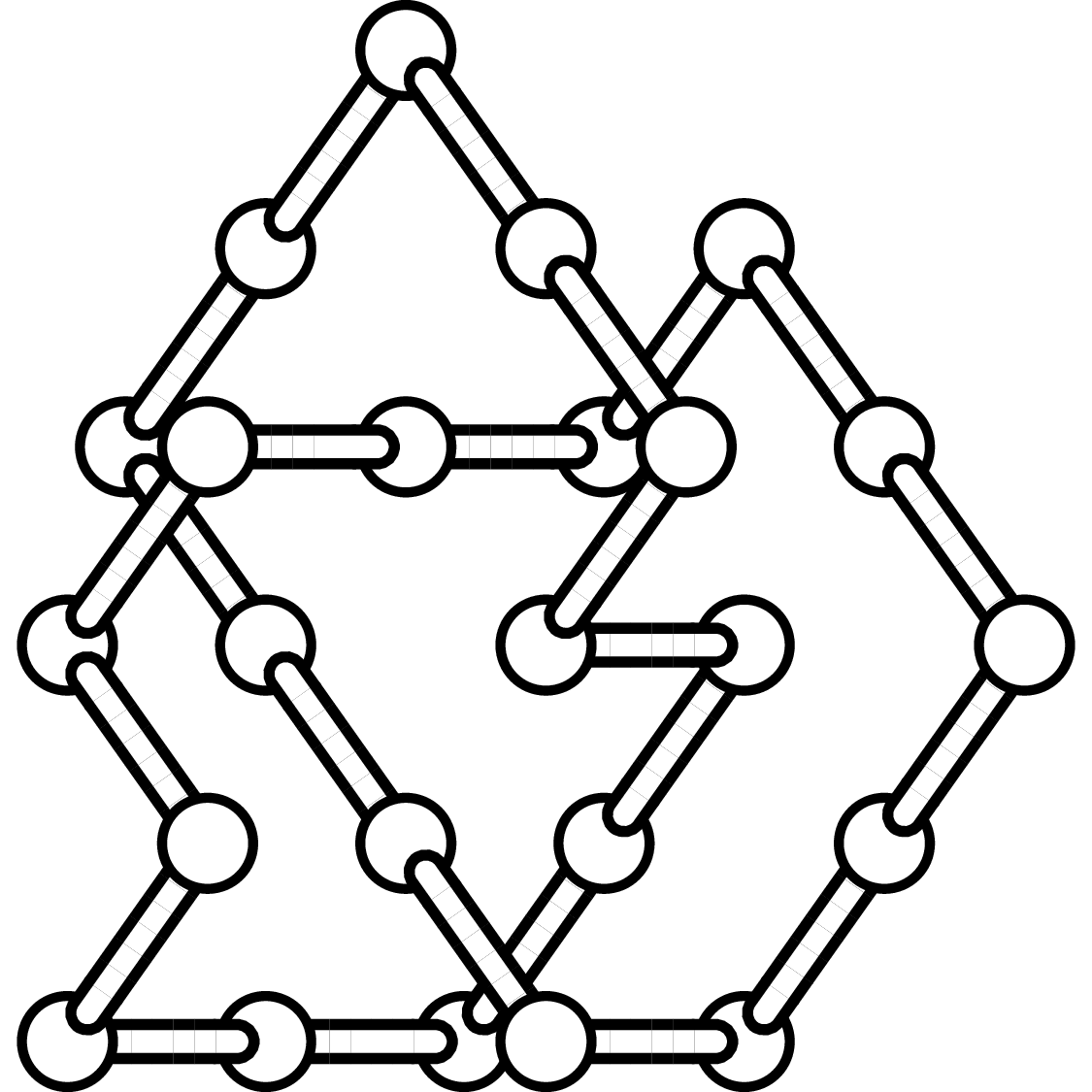}{Cubic lattice}\hss}
  \caption{The  trefoil  in various minimal edge forms.}
  \label{fig:EPSPL}
\end{figure}

\hypertarget{DiagramsSection}{}

\section{Making knot diagrams}

\label{sec:eps}

KnotPlot can generate lightweight Encapsulated PostScript (EPS) files using something like\\
\zcomm{psout filename}\\
and there are wide variety of options.
A simple example is shown in Figure {\ref{fig:EPS}}, produced with the commands:\\
\zcomm{mode s}\\
\zcomm{load 9.2.37}\\
\zcomm{psout knotLine}\\
\zcomm{psmode = 41     \% default psmode mode is 40}\\
\zcomm{psout knotRolfsen}\\
\zcomm{psmode = 45}\\
\zcomm{psout knotOther}

If \kp{} is in piecewise linear mode (\code{mode cb}), the EPS figures look
a little different. The following commands result in the images from Figure \ref{fig:EPSPL}.\\
\zcomm{mode cb; cradius = 0.16}\\
\zcomm{load ms/3.1; fitto avlength 3}\\
\zcomm{psout mstref}\\
\zcomm{load mseq/3.1}\\
\zcomm{fitto avlength 3; psout mseqtref}\\
\zcomm{load mscl/3.1; about x 45; about y 45; psout mscltref}\\

\begin{figure}
	\ppplen 0.333\textwidth     \hbox to \textwidth{\hss
	       	\ppp{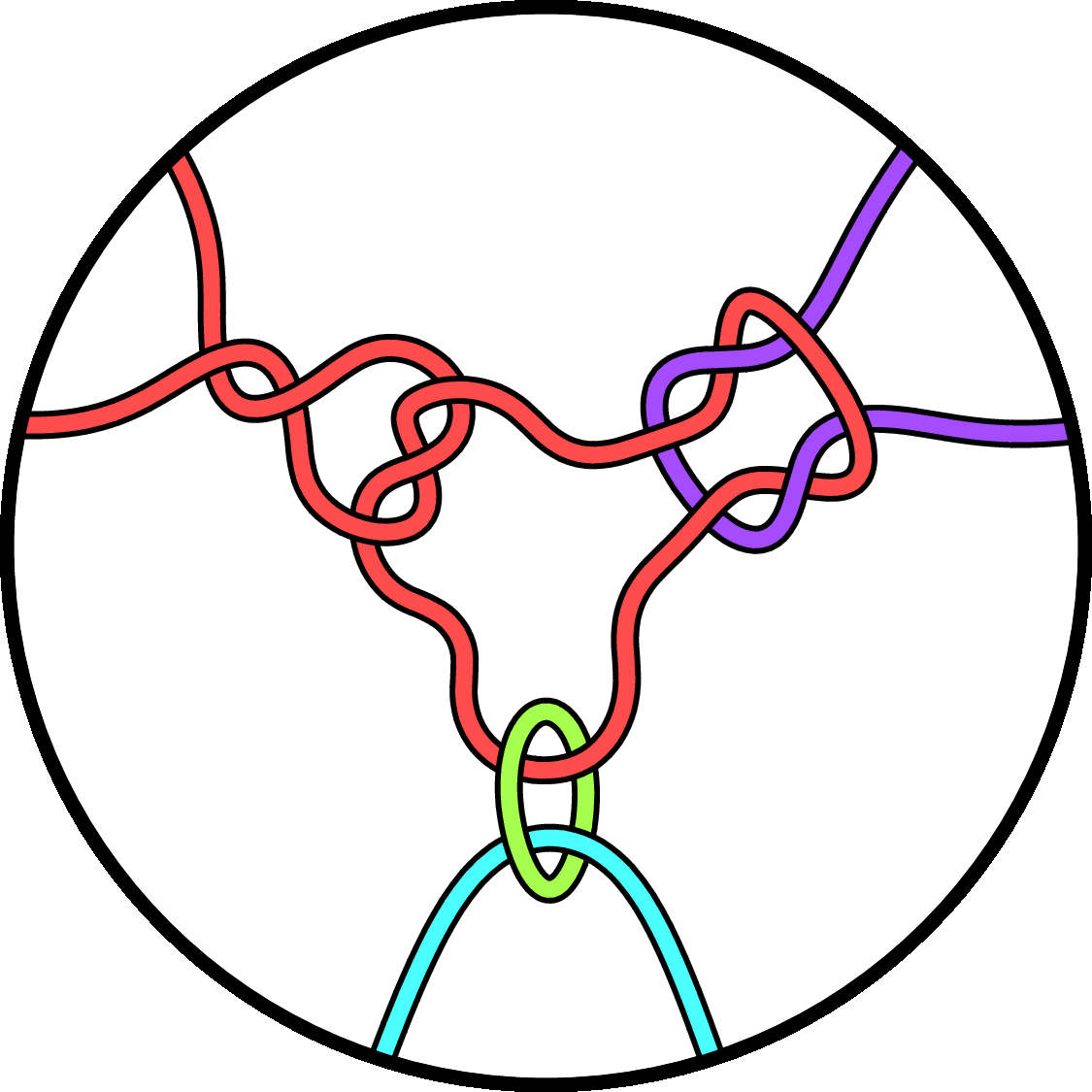}{three string tangle}\hss  
		\ppp{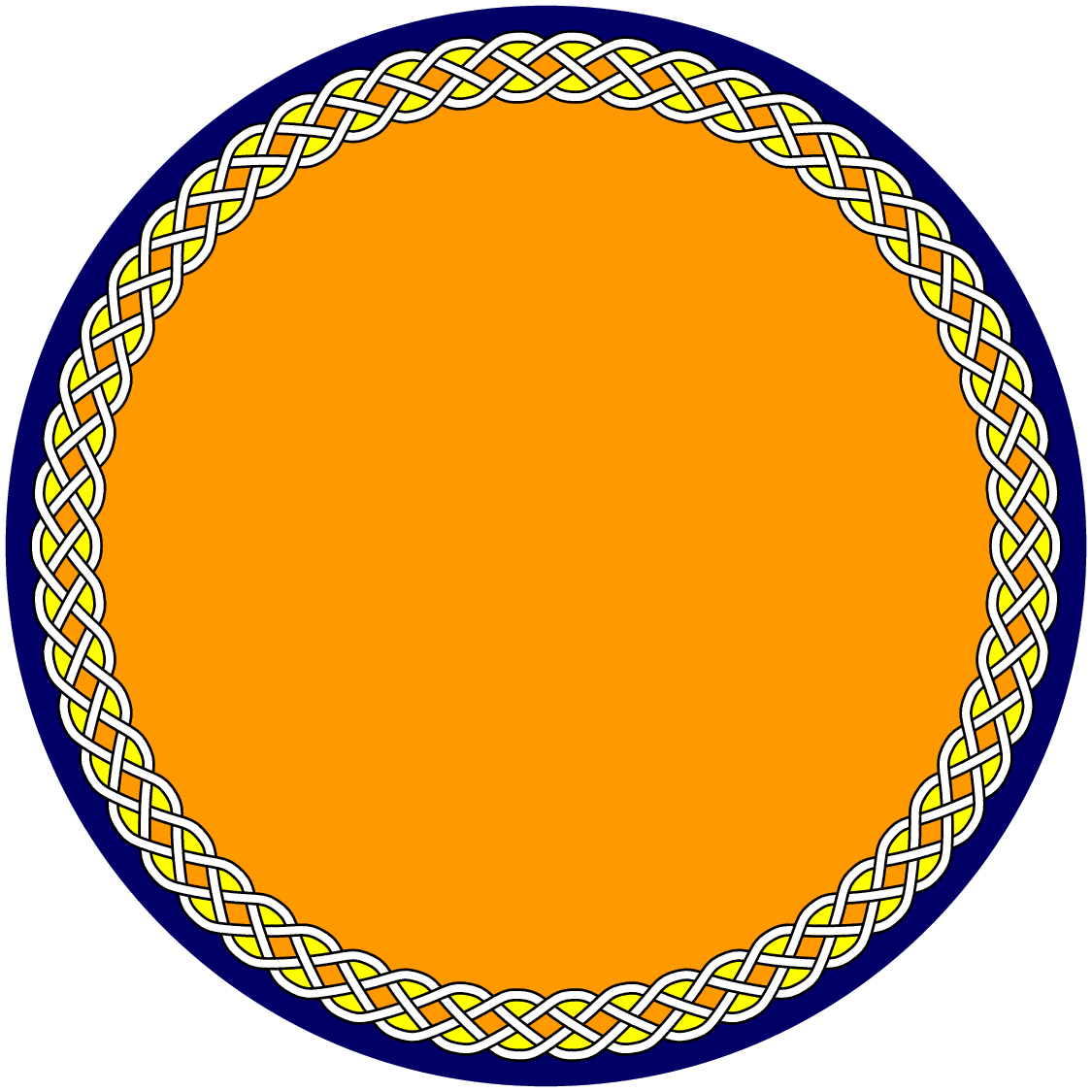}{Celtic knot}\hss
		\ppp{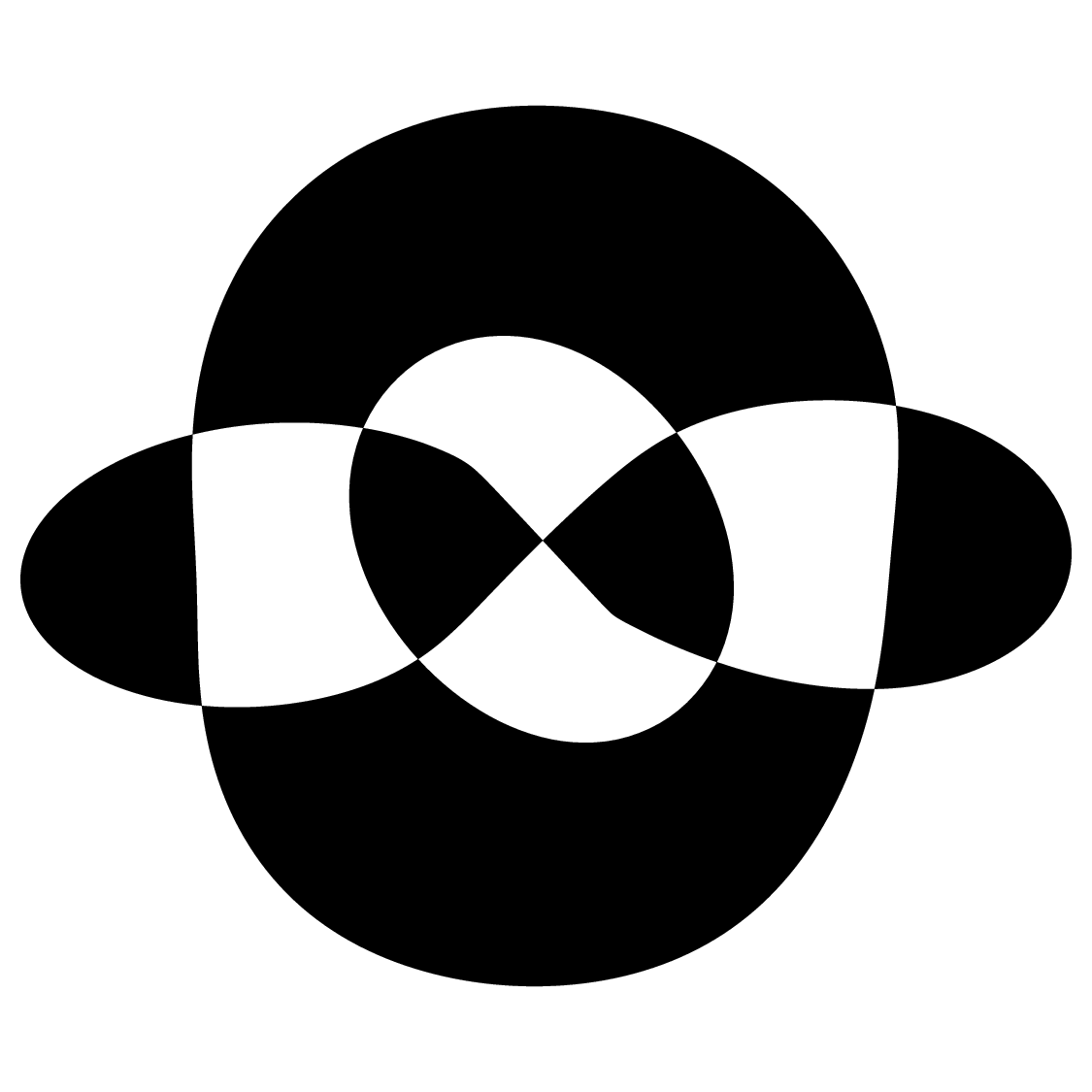}{checkerboard colouring}\hss} 
	\vspace*{3pt}
	\hbox to \textwidth{\hss
	       	\ppp{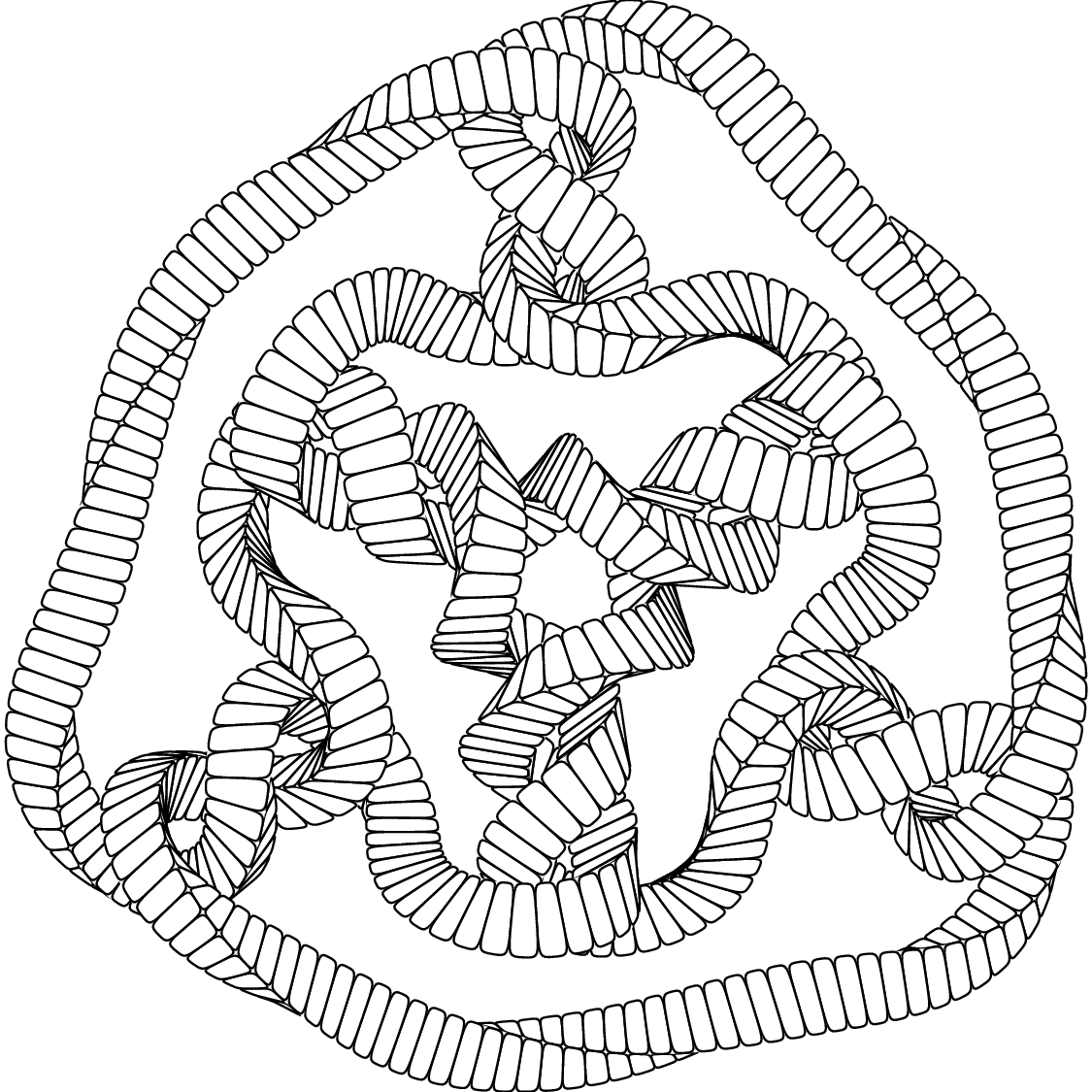}{Brunnian link}\hss  
		\ppp{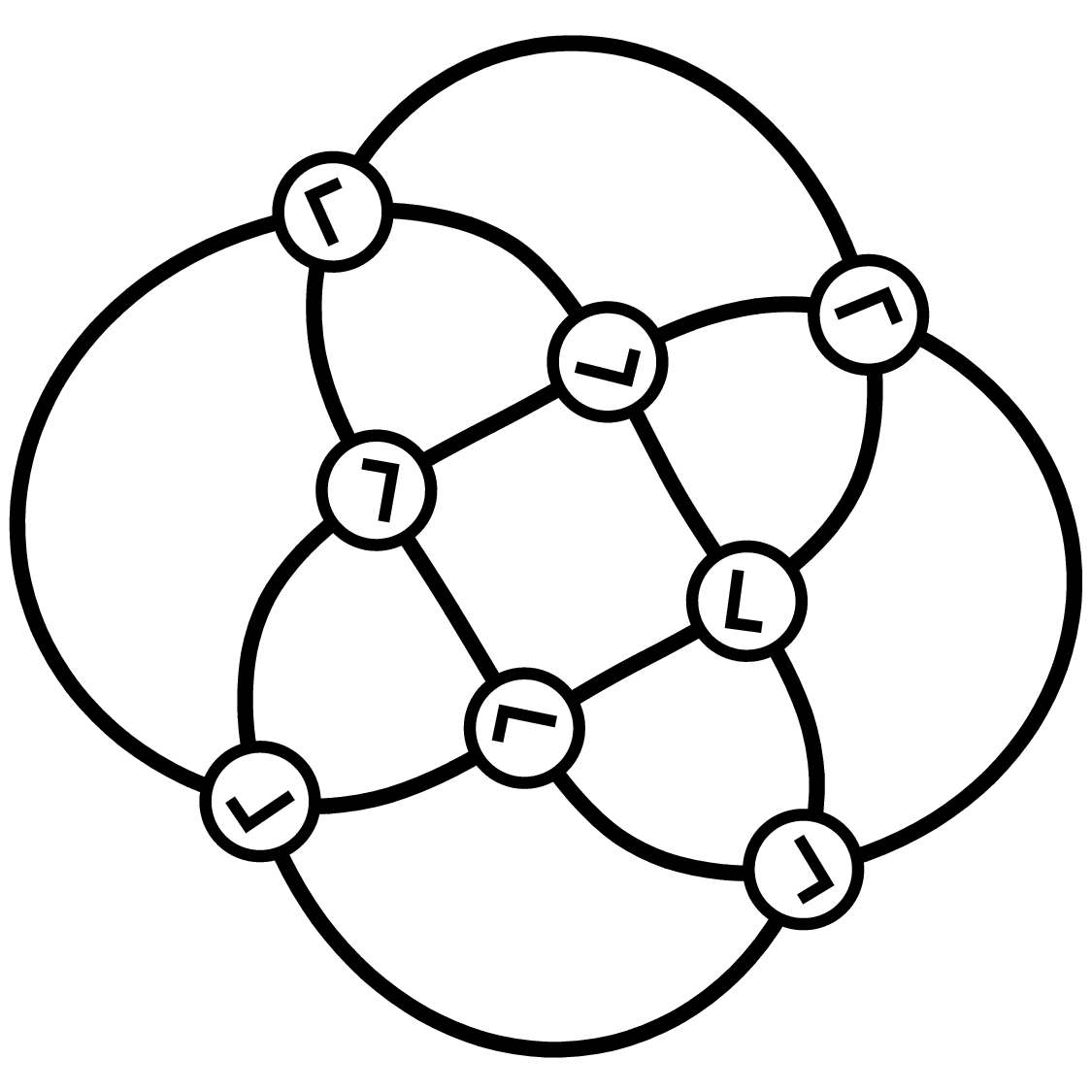}{tangle holder}\hss
		\ppp{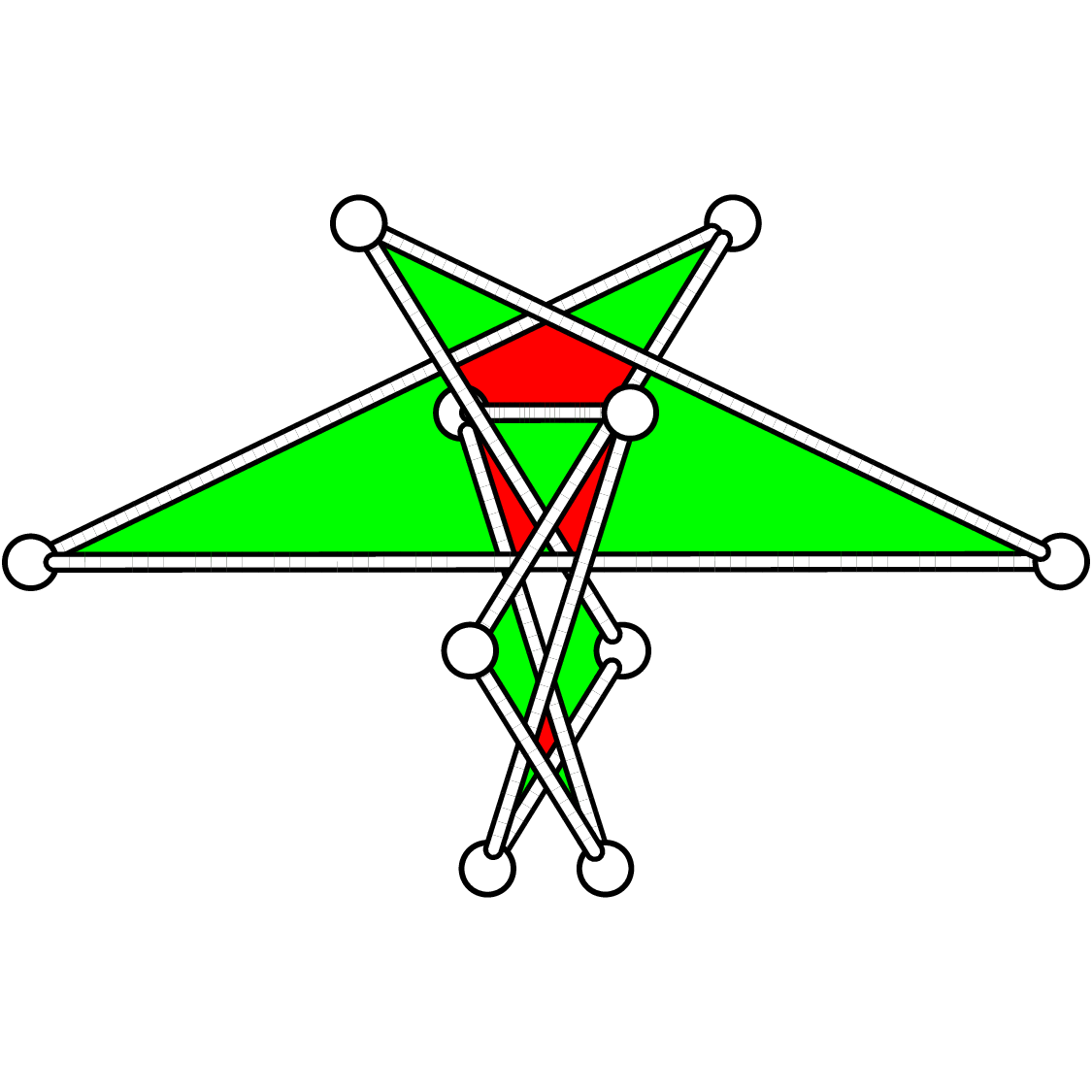}{minimal stick $10_{124}$}\hss} 

  \caption{Examples of Encapsulated PostScript diagrams}
  \label{fig:EPSexamples}
\end{figure}

RGS's dissertation \cite{Scharein_1998} was heavily illustrated with EPS images of this sort.
In addition to EPS images being lightweight, 
they also have the advantage as vector graphics that they may be scaled to arbitrarily large sizes without loosing their sharpness. 
To see the many flavours of PostScript output
together with the KnotPlot scripts that created them
visit the web page \cite{KnotPlotPostScript}.
A few of these are shown in Figure \ref{fig:EPSexamples} 
with the following links to the KnotPlot code that generated them:\\
\url{https://knotplot.com/postscript/m/TangleE.html}\\
\url{https://knotplot.com/postscript/m/CelticH.html}\\
\url{https://knotplot.com/postscript/m/psm_eofill1.html}\\
\url{https://knotplot.com/postscript/m/3compBrunA.html}\\
\url{https://knotplot.com/postscript/m/tangle-holder-B.html}\\
\url{https://knotplot.com/postscript/m/combine-eofill-PL.html}

\hypertarget{xconstab}{}
\hypertarget{creating}{}

\section{Creating knots and links}
\label{sec:creating}\def\exampleBraidWord{$  \relax (\sigma_{1}\sigma_{2}^{-1} )^{3}\sigma_{3}^{-1} \sigma_{1} \relax  $}
\def\exampleBraidWord{$  \relax (\sigma_{1}\sigma_{2}^{-1} )^{3}\sigma_{3}^{-1} \sigma_{1} \relax  $}
\def\exampleBraidWordNested{$  \relax \sigma_{5}(\sigma_{4}(\sigma_{2}\sigma_{3}^{-2})^{3}\sigma_{3}^{-1} \sigma_{1})^{2}\sigma_{2} \relax  $}
\def\exampleBraidWord{$  \relax (\sigma_{1}\sigma_{2}^{-1} )^{3}\sigma_{3}^{-1} \sigma_{1} \relax  $}
\def\exampleBraidWordNested{$  \relax \sigma_{5}(\sigma_{4}^{-1} (\sigma_{2}^{-1} \sigma_{3}^{2})^{3}\sigma_{3}^{-1} \sigma_{1})^{2}\sigma_{2} \relax  $}

\begin{figure}
	\hbox to \textwidth{\hss
	\qqq{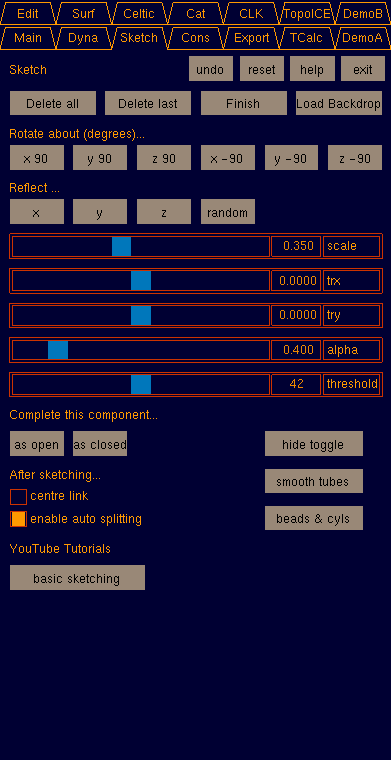}{\panel{Sketch} panel}\hss
	\qqq{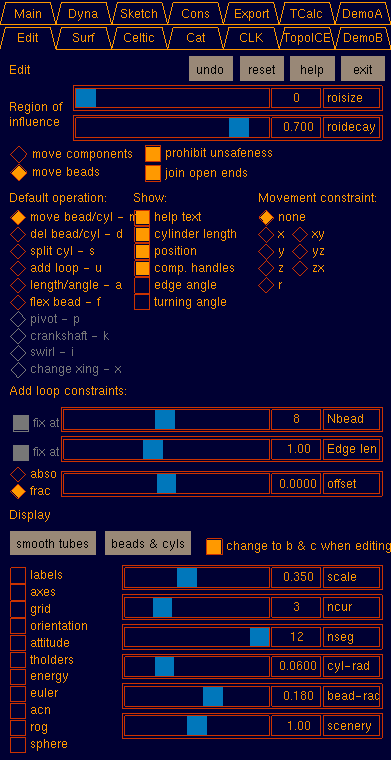}{\panel{Edit} panel}\hss
	\qqq{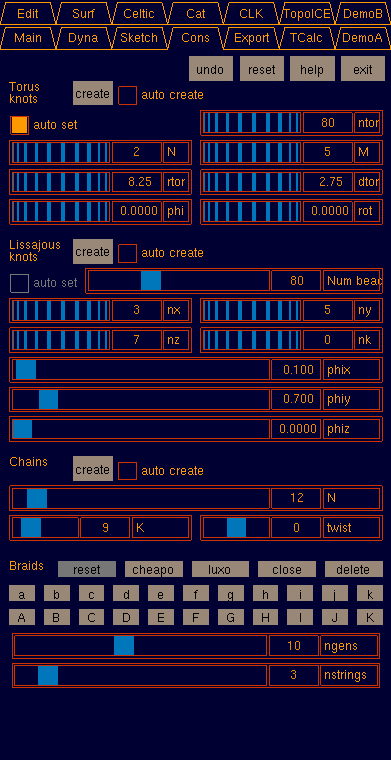}{\panel{Cons} panel}\hss}
\caption {Control panels for sketching, editing and constructing}
\label{fig:creating}
\end{figure}

\hypertarget{sketchtab}{}

\subsection{Sketching}
\label{sec:sketching}

If you have some knot you would like to have coordinates for, you can use
KnotPlot to sketch the entanglements.  If you click the  \panel{Sketch} panel, shown in Figure
\ref{fig:creating},  you can lay down points by left and right clicking.  When you
want to go over an existing strand, you use the right button.  To go
under, you use the left button.  You will see that as you lay points
down using the left and right buttons that the points will appear to
be further away and closer to you, respectively.  When you are done,
click either \pui{as open} or \pui{as closed} button to complete
the construction.  You can also just hold down a button and drag
around to lay down vertices at a reasonable pace.
Adjust the \pui{threshold} slider to control the spacing between beads 
when dragging.
When you are finished sketching, exit sketching by going to another panel, for example \panel{Main}.

When sketching, several keyboard operations can be useful, especially
when sketching on a trackpad.
\eric{fixing this\ldots}
\begin{itemize}
\item Hold down the SHIFT key when left clicking to simulate a right click that goes over.
\item Hold down the CONTROL key and left click and drag to translate the view.
\item Hold down the OPTION/ALT key and left click and drag to scale the view.
Click down at the location you wish to remain fixed when you scale.
\end{itemize}

\hypertarget{edittab}{}

\subsection{Editing}

If you want to move the vertices of a knot/link around, remove edges,
remove vertices, etc., use the \panel{Edit} panel.  There are lots of
choices here, but it is pretty straightforward to figure out.  If
\pui{prohibit unsafeness} is checked (the default), \kp{} will not let
you leave the knot in an unsafe position (see
\zoikt{advanceddynamics}{Section \ref{sec:dynamics}}) after any editing operation.
Note that if this box is checked, it is still possible to change the
knot type.  \kp{} does not check for unsafeness \emph{during} the
editing operation, only at the end.  To change the knot type, the
easier way is to leave the prohibit box checked, then increase
\pui{roisize} to a value near the top of the slider range.  Then click
and drag a bead or edge, pass a strand through another portion of the
knot, and when things look clear release the button.

If you hover the mouse button over a bead or an edge, some information
about that entity will be given.  Exactly what information depends on
what is selected in the \pui{Show:} check box section on the
\panel{Main} panel.  If you change
something on the \panel{Edit} panel, you might have to make sure the
Arena window has mouse focus before you seen any changes.  To do that,
just left click the mouse anywhere on the background and things should
work as expected.

The \pui{Default operation:} list in the \textit{Control Panel} shows
what editing operations are available and is initially set to
\pui{move beads/cylinders}.  If you want to delete one of those,
select \pui{del bead/cyl} and click on what you want gone.  Again,
remember to click on the background in the Arena to give that window
mouse focus.  It is actually not necessary to change the \pui{Default
operation:} to switch between operations.  You can press and hold one
of the shortcut characters \texttt{m}, \texttt{d}, \texttt{s},
\texttt{u}, \texttt{a}, or \texttt{f} to temporarily change the
operation (the corresponding shortcuts are shown to the right of the
operations).

You can constrain moving in various directions  by changing the
settings under the \pui{Movement constraint:} list on the right.
The letters indicate in which directions movement will not be allowed. 
For example, select \pui{xy} to change only the $z$-coordinates. 

One important thing to note is that, unlike when sketching, the
editing mode allows you to rotate, scale, or translate the view in the
normal \kp{} mouse movements.  Just remember when doing that, click down on an
empty space, avoiding the knot, and then drag the mouse.  It is for
this and other reasons that \panel{Sketch} panel will be removed in future
versions of \kp{} and its features added to the \panel{Edit} panel.

\def\dcomm#1{\hspace{0.4pc}\code{#1}\hspace*{0.4pc}}

Something fun to try:
\begin{enumerate}
\item Whilst on the \panel{Edit} panel, enter the command \dcomm{delete all}
\item Enter the command \dcomm{unknot}
\item Click on the \pui{beads \& cyls} button in the \pui{Display} section of \panel{Edit}
\item Enter command \dcomm{cut pieces 10}
\item Increase \pui{roisize} to at least 8
\item Click the \pui{none} diamond in the \pui{Movement constraint:} list
\item Click the \pui{join open ends} box until it is selected
\item Enter the command \dcomm{jitter 5}
\item Enter the command \dcomm{collision allow}
\item Enter the command \dcomm{until safe "ago 55"}
\end{enumerate}
You should now have a bunch of short strings that you can move around.
If you bring two open ends within a certain distance of each other, you should see a green
 indicator showing a connection between the ends.
Keep moving and the indicator will disappear.  However, if you release the mouse button when the indicator is on, the two beads will be connected
forming a new component. 
See if you can join all the ends to make a knot or a link.  
You might have to rotate and scale the view a few times. 
It might help to switch between \pui{move beads} and \pui{move components}.

\subsection{Torus knots and links}

The top part of the \panel{Cons} panel shown in Figure
\ref{fig:creating} is used to create torus knots and links.  To get a
quick start, check the \pui{auto create} box and start to adjust the
\pui{N} and \pui{M} rollers.  This will create an ($N$, $M$) torus
knot/link with \pui{ntor} edges, a meridional radius of \pui{rtor},
and a longitudinal radius of \pui{dtor}.  You can also leave \pui{auto
create} unchecked, adjust all the rollers to desired values, and then
click the \pui{create} button to create the knot.  As long as the
\pui{auto set} is checked, \kp{} will pick an appropriate value for
\pui{ntor}.  It can be interesting to uncheck this box, click the
\pui{beads \& cyls} button on \panel{Main}, and then set \pui{ntor} so
that the torus knot/link is undersampled.  You can use the
\hyperlink{torus}{\commandname{torus}} command to specify directly all
values.

\subsection{Lissajous knots}

Lissajous knots \cite{Lissajous} can be created with the next section
down on \panel{Cons} panel.  Things work in a similar manner to the
construction for torus knots, except the auto set is not available.
Note that many embeddings created in this construction will be in an
unsafe position (for a discussion of safety, see
\zoikt{advanceddynamics}{Section \ref{sec:dynamics}}).

The \demo{Lissajous} demo on \panel{DemoA} can be quite a lot of fun.
It chooses random values for the Lissajous knots and does not stop
until it finds a safe embedding.  When it does, you can click on the
\pui{simplify knot} button in the Arena.  Often what looks like a high
crossing number knot simplifies to an unknotted round circle.

\subsection{Chains}

This creates a chain with \pui{N} components with \pui{K} edges each
that twist.  More control can be had by using the
\hyperlink{chain}{\commandname{chain}} command.  To see what can be
done with chains, run the \demo{chains} demo on \panel{DemoA}.

\subsection{Braids}

In this section, you set the number of generators and strings and then
add positive or negative crossings one by one.  Once again, the use of
the \comm{braid} command is much more flexible.  To see what can be
done with braids, run the \demo{braids} demo on \panel{DemoA}.

\hypertarget{commands}{}

\section{Useful commands}
\label{sec:commands}

Below is a list of commands that are used with some frequency.
We have included example numbers in these, but all of the numbers
are specific to the context and what you are trying to do.
The current version of \kp{} has 357 commands,
please consult the manual \cite{KPM} to see what is not covered here.

\hypertarget{about}{\commandname{about x 45}} --- This rotates the knot/link while keeping
the axes fixed in place.  This changes the coordinates of
the knot/link.  The $x$ is the axis you are spinning about, so
you can use $y$ and $z$ too.  See also \hyperlink{rotate}{\commandname{rotate x 45}}.

\hypertarget{acn}{\commandname{acn}} --- Computes the exact average
crossing number  (up to computer precision), i.e., the average of the number of crossings seen over
all projections.  This also computes the writhe, which is the average over
all projections of the sum of the signed crossings.

\hypertarget{ago}{\commandname{ago 1000}} --- This runs the
dynamics/relaxation for 1000 steps.  While it is running, \kp{} will
not respond to any user actions.  The `a' in the command name stands
for `atomic', that is the command completes before anything else is
allowed to happen.  This should be used in place of
\hyperlink{go}{\commandname{go}} in scripts that are run in graphics
mode.  In non-graphics mode (see \zoikt{scripting}{Section \ref{sec:scripting}}) the
two commands are identical.  See also
\hyperlink{go}{\commandname{go}}.

\hypertarget{alias}{\commandname{alias plpsout "psout \$0"}} --- Creates an alias.  Aliases whose name begins with a $\sim$
are deleted upon a \comm{reset}.  Others, such as this one, will persist. 

\hypertarget{alignaxes}{\commandname{align axes}} --- This puts the
centre of mass at the origin and rotates the knot so that the
principal moment of inertia axes are aligned with the $x$-, $y$-, and
$z$-axes.  This command tends to produce embeddings that minimize
crossing numbers (in an orthographic projection to the $xy$-plane).

\hypertarget{allocate}{\commandname{allocate 10000}} --- This allocates   memory for  knots/links with larger numbers of beads.
Use with no argument to find the current amount allocated.

\hypertarget{angle}{\commandname{angle}} --- This computes the minimum
internal angle, maximum internal angle, ratio of maximum to minimum
internal angle, and the average internal angle for each component.

\hypertarget{angleturnng}{\commandname{angle turning}} --- 
Toggles between computing the internal angle and the turning angle. 
The turning angle, which measures
how much one edge turns relative to the one before it  is equal to 180 minus the internal angle.

\hypertarget{braid}{\commandname{braid}} \verb|(aB)^3Ca| --- Creates a braid with the braid word \exampleBraidWord.
The parentheses can be nested, so \code{braid} \verb|e(D(Bc^2)^3Ca)^2b| would produce the braid \exampleBraidWordNested.

\hypertarget{celticbox}{\commandname{celtic box 4 6 2}} --- Creates a Celtic box of the given dimensions.

\hypertarget{centre}{\commandname{centre}} --- Translates knot/link so that its axis-aligned rectangular bounding box  is at the origin.
The US spelling \commandname{center} is an alternative.

\hypertarget{centremass}{\commandname{centre mass}} --- Translates the knot/link so that its centre
of mass is at the origin. 

\hypertarget{chain}{\commandname{chain 10}} --- Creates a chain of unknots with 10 components.

\hypertarget{cheapo}{\commandname{cheapo}} --- This command
changes \param{ncur} and \param{nseg} to values lower than normally
used.  This command is useful when you want to reduce the rendering quality temporarily 
so that rotating and scaling the view is more responsive.
 The name borrowed from \cite[p.~150]{MetafontBook}.
There is also a complementary \comm{luxo} command.

\hypertarget{close}{\commandname{close 4}} --- Closes the indicated
component if the component is open.  Use \code{close all} to close all
components.  See the \comm{open} command.

\hypertarget{collision}{\commandname{collision allow}} --- Turns off KnotPlot's collision checking, allowing strand passages
during relaxation.  Use \code{collision fast} to go back to the default collision mode.

\hypertarget{colour}{\commandname{colour 0 red}} --- This colours
component 0 red.  Use \code{colour all red} to set
the colour to red for all components.  The available colours are in the
\code{resource/colours.txt} file within your \kp{} installation
directory.  The colour names are essentially those from X11
\cite{X11_color_names} with a few additions.  The names are case and
space insensitive.  You can also set the colour explicitly using RGB
values, for example \code{colour 2 rgb:.3/.5/1} where the values are
in the range 0 to 1.  For more complete information see
\cite{KPcolour}.  The US spelling \commandname{color} may be used as
an alternative.  See also \hyperlink{matrgbbead}{\commandname{matrgb
bead}} and \hyperlink{matrgbknot}{\commandname{matrgb knot}}.

\hypertarget{conway}{\commandname{conway 6,3,4}} --- Creates a pretzel
link from the Conway notation \cite{Conway}.  Also prints out the
\comm{tangle} command string used to generate the link in the tangle
calculator, in this case \verb|tangle 6r3r#4r#N|

\hypertarget{cut}{\commandname{cut 27}} --- Deletes the edge starting
at a given bead/vertex.  Use \code{show labels} or check the box on
\panel{Main} to see the bead labels. 
Vertices are numbered starting at 0 and
edge $n$ goes from vertex $n$ to vertex $(n + 1)$ (modulo the number
of edges, if the configuration is closed).

\hypertarget{cutoutside}{\commandname{cut outside x 0}} --- Cut all edges that are entirely 
\eric{outside a plane is standard computer graphics terminology.  a plane has its normal, outside means on the side of the plane 
where in the direction of the normal}
outside the plane $x=0$.
To see all the options, run the \demo{cut} demo in the Lua scripted demos section of \panel{DemoA} (also view the Lua source).

\hypertarget{cwd}{\commandname{cwd}} --- This tells you  the current working  directory that \kp{} is running
from.  If you save something, it goes into this
directory.  See also the  \hyperlink{path}{\commandname{path}} command.

\hypertarget{delete}{\commandname{delete} \code{1}} --- If you have a link, this will delete
component 1 (components are numbered starting at 0).
Use \comm{undo} if that was a mistake. 
See also \hyperlink{keep}{\commandname{keep 1}}.

\hypertarget{deletedownto}{\commandname{delete downto 15}} --- Deletes vertices, while keeping the
knot type fixed, to 15 vertices.  Oftentimes it will not be able to
get to that value, in which case it just stops.  You can say
\commandname{delete downto 0} 
(equivalent to the \pui{dbeads} button on \panel{DemoA})
and that will do the best it can to get
rid of extra vertices.  
Doing this periodically during relaxation and with the parameter \param{stusplit} set to a positive value
is how the knots on the \demo{Split/delete} demo on \panel{DemoA} work.
See the exercise in Section \hyperlink{simple}{\ref{sec:simple}}.

\hypertarget{diagram}{\commandname{diagram 8 -10 -16 12 -14 6 2 -4}} --- Generates a configuration whose
Dowker code is 8 -10 -16 12 -14 6 2 -4. 
Uses   (with permission) the knot diagram engine
 by Morwen Thistlethwaite and Kenneth Stephenson modified slightly from the `draw' routine 
  in Knotscape \cite{knotscape}.

\hypertarget{dispt}{\commandname{display true}} --- Forces a redisplay of the current state of the KnotPlot Arena.
Useful mostly in scripts run in graphics mode when you want to ensure the display is correctly updated before creating a picture.

\hypertarget{distance}{\commandname{distance}} --- Shows the minimum
distance between non-adjacent edges.  This is the value reviewed when
you check on safety with the \comm{safe} command.  See 
\zoikt{advanceddynamics}{Section \ref{sec:dynamics}}.

\hypertarget{dowker}{\commandname{dowker}} --- KnotPlot computes
Dowker codes \cite{DowThis} with this command as well as Extended Gauss
Codes (EGCs) \cite{ghomfly} using the \comm{gauss} command.  By
default, these commands compute their codes in a fixed projection that
may not correspond to the view in the Arena.  Please see the
discussion of the \hyperlink{dowkerprojectionz}{\commandname{dowker
projection}} options on how to change modes so that the codes match the
view.  For a comprehensive introduction to crossing codes see the
\demo{codes} demo on \panel{DemoA}.


\hypertarget{dowkerprojectionview}{\commandname{dowker projection
view}} --- Set the Dowker projection mode to the view (as
an orthographic projection) of the knot/link seen in the \textit{KnotPlot
Arena}.  In this
mode, the \hyperlink{dowker}{\commandname{dowker}} and
\hyperlink{gauss}{\commandname{gauss}} commands will correspond to
what is seen in the view window, and will change if you rotate the
view.
Same as \hyperlink{gaussprojectionview}{\commandname{gauss projection view}}.

\hypertarget{dowkerprojectionz}{\commandname{dowker projection z}} ---
Set the Dowker projection mode to the $z$-projection (i.e., onto
the $xy$-plane).  In this mode,
crossings are computing by orthographically projecting the edges to
the $xy$ plane.  The \hyperlink{dowker}{\commandname{dowker}} and
\hyperlink{gauss}{\commandname{gauss}} commands will give codes that,
if the knot has been rotated, do not correspond to what you see in the view window,
especially if the knot is rotated.  Alternately, one can do
a \hyperlink{rotatefix}{\commandname{rotate fix}} followed
by \hyperlink{dowker}{\commandname{dowker}}.
Same as \hyperlink{gaussprojectionz}{\commandname{gauss projection z}}.

\hypertarget{drawflat}{\commandname{draw flat}} -- Draws the knot/link so that the cylinders
look like they are made out of rectangles.

\hypertarget{drawhflat}{\commandname{draw hflat}} -- Draws the knot/link so that the cylinders
look like they are made out of long rectangles.

\hypertarget{drawnormal}{\commandname{draw normal}} -- Draws the knot/link in the default fashion.

\hypertarget{drawspectrum}{\commandname{draw spectrum}} -- Draws the knot/link with rainbow colouring.

\eric{added an example}

\hypertarget{duplicate}{\commandname{duplicate}} --- Duplicates component 0.  The duplicate will be directly on top of the original so you
will not see both unless you move one of them.
An example usage is to quickly create a composition of two $6.3$ knots:\\
\zcomm{load 6.3}\\
\zcomm{shift maxx  \% shifts bead 0 to maximum x-position}\\
\zcomm{duplicate}\\
\zcomm{open all}\\
\zcomm{reflect x 1; translate 20 0 0 1}\\
\zcomm{join F0 F1; close all}\\
\zcomm{ago 333; centre}\\

\hypertarget{energy}{\commandname{energy}} --- Prints the knot energy according to the current energy model.

\hypertarget{energymodel}{\commandname{energy model}} --- List the available energy models (there are more than a dozen to choose from). 
\eric{list em all added}
Change the energy model used with commands such as
\begin{itemize}
\item \commandname{energy model MD} --- Simon's minimum distance energy \cite{SimonMD}
\item \commandname{energy model symm} --- Buck and Orloff's symmetric energy \cite{BuckOrloffB}
\item \commandname{energy model nbeads} --- Energy is the total number of beads in the knot. Can be useful when trying to reduce the number of vertices.
\end{itemize}
Consult the KnotPlot Manual for more information~\cite{KPM} on the different models.

\hypertarget{exit}{\commandname{exit}} --- Same as the \comm{quit} command.
This command will quit \kp{} if run from the \textit{Command Window}.
In scripts this causes the script to quit at that point in the script.  

\hypertarget{export}{\commandname{export mymodel}} ---  Exports a \3\ polygonal mesh to the file \code{mymodel.obj}
which is in Wavefront OBJ format~\cite{OBJ}.

\hypertarget{fitto}{\commandname{fitto 10}} ---  Multiplies the coordinates of the knot and centres it at the origin such that
the maximum extent in the $x$, $y$ or $z$ directions is $\pm10$.  This command
is useful for rescaling and centering the model so that it is in the viewing window.

\hypertarget{fittomindist}{\commandname{fitto mindist 0.15}} --- This is like the normal
\hyperlink{fitto}{\commandname{fitto}} but it makes it so that the knot/link has a
minimum distance of  0.15 between any two non-adjacent cylinders.


\hypertarget{gauss}{\commandname{gauss}} --- Works the same as the \hyperlink{dowker}{\commandname{dowker}} command
except that it computes the extended gauss code.
How EGCs are computed is explained in Section~2.3
of the HOMFLY-PT chapter in this book. 

\hypertarget{gaussprojectionview}{\commandname{gauss projection view}} --- 
Same as \hyperlink{dowkerprojectionview}{\commandname{dowker projection view}}.

\hypertarget{gaussprojectionz}{\commandname{gauss projection z}} --- 
Same as \hyperlink{dowkerprojectionz}{\commandname{dowker projection z}}.

\hypertarget{go}{\commandname{go 1000}} --- This runs the dynamics/relaxation for
1000 steps.  
While it is running, \kp{} will respond 
commands and user interactions. 
See also \hyperlink{ago}{\commandname{ago}}.

\hypertarget{head}{\commandname{head 22}} --- Shows only the first 22
beads of the knot and hides the rest.  Used on the demo \demo{braid
warp} on \panel{DemoA}.  Although the beads are hidden, they are
really still there.  Enter \code{head off} turn off hiding.

\hypertarget{hide}{\commandname{hide 1}} --- Hides component number 1.
It is still there, but you will not be able to see it.  The opposite
is \hyperlink{unhide}{\commandname{unhide 1}}.

\hypertarget{hideall}{\commandname{hide all}} --- Hides all of the
components.  They are still there, but you will not be able to see
them.  The opposite is \hyperlink{unhideall}{\commandname{unhide
all}}.

\hypertarget{homfly-pt}{\commandname{homfly-pt}} -- Computes the HOMFLY-PT
polynomial.  
Alternative command names are \commandname{homfly}, \commandname{homflypt} and \commandname{flypmoth}.
For references and a description of a stand-alone tool to compute these polynomials
please refer to the HOMFLY-PT chapter in this book.

\hypertarget{id}{\commandname{id}} -- Attempts to identify the currently loaded knot or link using
the \hyperlink{homfly-pt}{\commandname{homfly-pt}} command  and KnotPlot's built-in lookup table for
HOMFLY-PT polynomials. 
This lookup table is limited in scope and not changeable by users.  The HOMFLY-PT chapter 
has links to various knotting tools that have much more flexibility. 

\hypertarget{imgout}{\commandname{imgout tref}}  --- Saves the contents of the view window as the PNG file \code{tref.png}.
If you want a JPEG instead, append \code{.jpg} to the file name, e.g.,
use the command \code{imgout tref.jpg}. 

\hypertarget{info}{\commandname{info}} --- This tells the number of components, number of
vertices/beads, and some other information.
Use \code{info s} to get a short version of the same information.

\hypertarget{jitter}{\commandname{jitter 0.1}} --- This moves each vertex by a random value
$\leq 0.1$.  If you have a lot of parallel edges (like in a lattice
knot) and do something like try to find a crossing code, you will
run into all sorts of problems (you will get division by 0 in some
of the computations).  You can use \commandname{jitter} to get
out of that problem.

\hypertarget{join}{\commandname{join L0 F1}} --- Join the last bead of component 0 to the first bead of component 1, creating
a new component.   Both components must be open.

\hypertarget{keep}{\commandname{keep 1}} --- If you have a link, this will delete all of
the components but component 1.  
See also \hyperlink{delete}{\commandname{delete 1}}.

\hypertarget{knot}{\commandname{knot}} --- Randomly choose a knot/link
from the Rolfsen catalogue and load it.  To choose a random link with 
3 components use \code{knot 3}.

\hypertarget{length}{\commandname{length}} --- This computes the total
arc length, minimum edge length, maximum edge length, ratio of maximum
to minimum edge lengths, and the average edge length for each
component.

\hypertarget{line}{\commandname{line from 5 3 1 to -4 2 -7}} --- Create a line between two locations.  You can omit the \code{from}
and \code{to}. 

\hypertarget{lnknum}{\commandname{lnknum}} --- If you have a link, this will
compute the linking number matrix, i.e., the linking number between each pair
of   components.

\hypertarget{load}{\commandname{load file}} --- Loads a knot/link file.
See Section \hyperlink{LoadNSave}{\ref{sec:loadsave}}
for information  on using this command.

\hypertarget{loadcombine}{\commandname{load combine file}} --- If you have one knot/link
loaded and you want to load another knot/link as to form a union with what is
there, you use \commandname{load combine}.

\hypertarget{loadsum}{\commandname{load sum file}} --- This is used
if you have a knot/link loaded
and you want to do a composition with a new knot/link.  You usually
need to relax the knot/link after doing this.  For links, you can
specify which component of the currently loaded link that you
connect sum with.  For example, if you\\
\hyperlink{load}{\command{load 9.3.2}}\\
\command{load sum 3.1 0}\\
you get something different than if you\\
\hyperlink{load}{\command{load 9.3.2}}\\
\command{load sum 3.1 1}\\
For links, the component to compose with is set to a default value of 0.

\hypertarget{lua}{\commandname{lua run Brunnian}} --- Run the given Lua program, in this case
the \demo{Brunnian} demo on \panel{DemoA}.

\hypertarget{luxo}{\commandname{luxo}} --- This command
changes \param{ncur} and \param{nseg} to values higher than normally
used.  This command is useful before exporting an image with
\comm{imgout}.  The name borrowed from \cite[p.~150]{MetafontBook}.
There is also a complementary \comm{cheapo} command.
\eric{cheapo added}

\hypertarget{massopen}{\commandname{mass open}} --- If you have an open knot, this anchors the
first and last vertex so that they do not move during relaxation.  

\hypertarget{matrgbbead}{\commandname{matrgb bead 0.1 0.4 0.8}} --- This changes the beads to the
colour with RGB values of 0.1, 0.4, and 0.8.  This affects all
components.

\hypertarget{matrgbknot}{\commandname{matrgb knot 0.5 0.25 0.3}} --- This changes the cylinders
to the colour with RGB values of 0.5, 0.25, and 0.3.  This affects all
components.  If you just want to change the colour of one component,
see \hyperlink{colour}{\commandname{colour}}.

\hypertarget{mf}{\commandname{mf myknot}} --- Outputs the current knot 
as a Metafont \cite{MetafontBook} program.  Only works for 1-component knots.
\eric{no, not for links}

\hypertarget{modecb}{\commandname{mode cb}} --- Puts \kp{} in \commandname{beads and cylinders} mode.
There is a button for this on the main panel, so you would probably only
use this in a script.

\hypertarget{modes}{\commandname{mode s}} --- Puts \kp{} in \commandname{smooth tubes} mode. There is a
button for this on the main panel, so you would probably only use this in a
script.

\hypertarget{nap}{\commandname{nap} $N$} --- Pauses for $N$ seconds.  You might use this
if you are running a script to view knots, kind of like a slideshow,
and you want it to pause a bit between each knot.

\hypertarget{nbeads}{\commandname{nbeads +5}} --- Adds 5 to the number of vertices.  It
does so by following a spline of the knot, so the new version
typically has nicely distributed edge lengths.  You can also subtract
using \commandname{nbeads -7}, or specify a number exactly using
\commandname{nbeads 20}.  The latter is precarious because you can
easily change knot type.  You really only want to use this
command if your knot/link is pretty smoothish looking to start.
See also \hyperlink{split}{\commandname{split}} and \hyperlink{nbeadsmult}{\commandname{nbeads mult}}.

\hypertarget{nbeadsmult}{\commandname{nbeads mult 3}} --- This is like \hyperlink{nbeads}{\commandname{nbeads}} but
works as a multiplier.  So this example triples the number of
vertices.  You really only want to use this
command if your knot/link is pretty smoothish looking to start. See also \hyperlink{nbeads}{\commandname{nbeads}} and
\hyperlink{split}{\commandname{split}}.

\hypertarget{open}{\commandname{open 3}} --- Opens  component 3.  
Use \code{open all} to open all components.   See the \comm{close} command.
This command opens a closed component by removing the edge between the 0$^\text{th}$  bead and the
 last bead.  If you want to open the component at another
edge, see the \hyperlink{shift}{\commandname{shift}} command.

\hypertarget{orthographic}{\commandname{orthographic}} --- Sets the
view to an orthographic projection.  In the default
\hyperlink{perspective}{\commandname{perspective}} view, objects that
are further away look smaller.  In the orthographic view, this is not
the case.

\hypertarget{panel}{\commandname{panel set Cons}} ---  Set the Control Panel to the indicated tab.  This could be used in an initialization script, for example,
to get \kp{} into a desired setup.

\hypertarget{parameters}{\commandname{parameters}} --- Lists all of the parameter values.
If an argument is given, then \kp{} will list all the parameters that start with that prefix.
For example \code{parameters ps} will list all the parameters dealing with EPS output.
You can do a\\
\command{parameters > pars.txt}\\
and then open the \texttt{pars.txt} file to see all of the parameters.

\hypertarget{path}{\commandname{path}} --- Prints the read, write, and execute paths.  See Section \hyperlink{aloadsave}{\ref{sec:aloadsave}}.

\hypertarget{perspective}{\commandname{perspective}} --- This switches
the view to a perspective view, which is the default.  In the
projection view, things that are further away look smaller.  In the
orthographic view, this is not the case.  The alternative is
\hyperlink{orthographic}{\commandname{orthographic}}.

\hypertarget{project}{\commandname{project random}} --- Project the knot in a random direction. 
Changes the embedding, not the view.

\hypertarget{psoption}{\commandname{psoption bbox square}} --- Sets the bounding box option for \comm{psout}.
The are many other options that can be set.
See the KnotPlot code on the PostScript examples page \cite{KnotPlotPostScript}.

\hypertarget{psout}{\commandname{psout diagram}} --- Generates an eps
file called \code{diagram.eps} of a knot diagram version of the
configuration you have loaded.  See
\hyperlink{DiagramsSection}{Section \ref{sec:eps}}.

\hypertarget{quit}{\commandname{quit}} --- The same as the \comm{exit} command.

\hypertarget{refine}{\commandname{refine 3}} --- Refine the knot using B\'ezier interpolation by a given factor
 ($3$ times as many beads in this case).
  This command may change the knot type.
  
 \hypertarget{refineequi}{\commandname{refine equilateral 1.02}} --- Refines the knot so it is approximately 
 equilateral with the specified edge length.  Use the \comm{length} command to see how well it refined.

\hypertarget{reflect}{\commandname{reflect x}} --- This reflects the knot/link relative
to the chosen $x$-axis.  You can also use $y$ or $z$.
In other words, it flips the signs of the
$x$ (or $y$ or $z$) coordinates.  Any one of these will create the
mirror image of the given knot/link. 
You can combine reflections, as in \code{reflect yx} or reflect in a randomly chosen direction \code{reflect r}.
If you want to reflect a single component of a link, specify 
the component number \code{reflect z 2}.

\hypertarget{reset}{\commandname{reset}} --- This resets the viewing transform, the display mode, and
other things,
keeping the knot/link currently loaded.

\hypertarget{resetall}{\commandname{reset all}} --- Resets \kp{} to its startup state. 
Also see the gentler \hyperlink{reset}{\commandname{reset}} command.

\hypertarget{revbeads}{\commandname{revbeads 1}} --- This reverses the orientation on
component 1.  
Without any arguments
\commandname{revbeads} reverses all the components.
 If you click on the \button{Labels} box on the right side of \panel{DemoA}, you can
see the vertex numbering change when you do this.

\hypertarget{rog}{\commandname{rog}} --- Computes the radius of gyration.

\hypertarget{rotate}{\commandname{rotate x 45}} --- This rotates the camera while keeping the
knot/link the axes fixed in place.  This does not change the
coordinates of the knot/link.  You need to use \hyperlink{rotatefix}{\commandname{rotate
  fix}} after this, or after rotating the knot/link with the mouse, if
you want to change the coordinates.  The $x$ is the axis you are
spinning about, so you can use $y$ and $z$ too.  See also
\hyperlink{about}{\commandname{about x 45}}.

\hypertarget{rotatefix}{\commandname{rotate fix}} --- If you have rotated the view using
\hyperlink{rotate}{\commandname{rotate}} or using the mouse, the coordinates of the
actual knot/link will not have changed, just the viewing transformation.
This command changes the coordinates so that the view you are seeing is
lined up so that you are looking down the $z$-axis.  
A useful command if have a nice view of something you created and want to adjust the
embedding so that when you save your knot and load it back into KnotPlot, you 
do not have to any rotating to get back to that view.

\hypertarget{rotateunit}{\commandname{rotate unit}} --- Resets any rotations you have done
with the mouse or with the \comm{rotate} command. 

\hypertarget{safe}{\commandname{safe}} --- This tells you whether
the knot/link is in a \kp{} safe configuration. 
See Section \hyperlink{advanceddynamics}{\ref{sec:dynamics}} for details.

\hypertarget{save}{\commandname{save file}} --- Saves a knot/link file.
See Section \hyperlink{LoadNSave}{\ref{sec:loadsave}}
for information  on using this command.

\hypertarget{scale}{\commandname{scale 2}} --- This multiplies all of the vertex
coordinates by 2.

\hypertarget{scalenonprop}{\commandname{scale 2 4 0.8}} --- This non-proportionally scales the knot/link coordinates
by a factor of 2 in the $x$-direction, 4 in the $y$-direction and 0.8 in the $z$-direction.  This command will not change the knot type.

\hypertarget{shift}{\commandname{shift 5}} --- This shifts the labels on the vertices by
5.  So, for example, if you want some other vertex to be the
$0^{th}$ vertex, you can use a shift to make that happen. 

\hypertarget{show}{\commandname{show orientation}} --- Shows
orientation arrows on the components.  This command has many other
options in addition to \code{orientation} (see the \pui{show} check boxes
on the \panel{Main} panel).

\hypertarget{split}{\commandname{split}} --- This adds a vertex half way into each edge, doubling
the number of vertices.  See also \hyperlink{nbeads}{\commandname{nbeads}} and \hyperlink{nbeadsmult}{\commandname{nbeads mult}}.

\hypertarget{swap}{\commandname{swap 2 4}} --- Swaps the numbering of the components indicated.
Use \code{swap random} to swap all components randomly.


\hypertarget{tangle}{\commandname{tangle}} \verb|21*r2r#2r#N| --- Create the link $7^2_5$ using the tangle calculator. 
Use the \comm{conway} command to see how the Conway notation gets parsed into tangle calculator commands.
Also try the interactive tangle calculator on the \panel{TCalc} panel.

\hypertarget{tfunction}{\commandname{tfunction 2.3 knot}} --- Queue an arbitrary command to be executed in a given number of seconds.

\hypertarget{thickness}{\commandname{thickness}} --- This computes some
sense of thickness.  
It is not equivalent to EJR's definition of knot thickness \cite{dissertation,meideal,mecancomputers}.

\hypertarget{timer}{\commandname{timer start gumby}} --- Start a named timer.  The name is arbitrary and
you can have as many timers as you wish.   Check the timer at some point in the future with \code{timer check gumby}.

\hypertarget{torus}{\commandname{torus 3 5 50 11 2}} --- Creates a (3,5) torus knot with
50 edges where the  longitudinal radius is 11 and the meridional radius is 2.  
You can skip last three parameters and \kp{} defaults to the current value of the parameters \texttt{N-torus},
\texttt{R-torus}, and \texttt{d-torus}.

\hypertarget{translate}{\commandname{translate 5 8 3}} --- This translates the knot 
by 5 in the $x$ direction, 8 in the $y$-direction, and 3 in the $z$-direction.  
Changes the embedding, not the viewing transform in a similar manner to   \hyperlink{about}{\commandname{about}}.
To un-transform a translation, do a \comm{undo} command immediately after the unwanted transformation.
Note, translations \emph{are not undone} by  \hyperlink{untran}{\commandname{untran}}, which resets the viewing transform.

\eric{clarified...}

\hypertarget{translateto}{\commandname{translate to 5 3 7}} --- This
translates the knot so that the 0th vertex is at (5,3,7).  You
can use \hyperlink{shift}{\commandname{shift}} to set a different
vertex as the 0th vertex if you want to move a different vertex
to a particular point in space.

\hypertarget{twfix}{\commandname{twfix}} ---  Sets the value of \code{twist} (found on the \panel{Surf} panel) to the
value required so that the cross sections of the smooth tubes align at the beginning and end.
For example, to see what this does, type\\
\hyperlink{load}{\command{load 8.18}}\\
\hyperlink{drawflat}{\command{draw flat}}\\
{\command{\param{nseg} = 3}}\\
You will find an unsightly twist in the smooth tube.
Entering \hyperlink{twfix}{\commandname{twfix}} makes this disappear.

\hypertarget{undo}{\commandname{undo}} --- This undoes any changes just made to the knot.

\hypertarget{unhide}{\commandname{unhide 1}} --- Unhides component number 1 if it is
hidden.  Now you can see it.  The opposite is \hyperlink{hide}{\commandname{hide 1}}.

\hypertarget{unhideall}{\commandname{unhide all}} --- Unhides all of the components.  Now you
can see all of the components.  The opposite is \hyperlink{hideall}{\commandname{hide
  all}}.

\hypertarget{unknot}{\commandname{unknot 50 5}} --- Creates a regular 50-gon on a circle
of radius 5.  
If the radius is omitted, it defaults to the current value of \param{N-torus}.

\hypertarget{untran}{\commandname{untran}} ---  Resets the viewing transform to the unrotated/untranslated one.

\hypertarget{version}{\commandname{version}} --- Use this command to find out the version, build number, 
and compile date of the KnotPlot you are running.
Also shows you where the KnotPlot home directory is.

\hypertarget{xing}{\commandname{xing}} --- This tells you the number of crossings seen
looking down the $z$-axis.  If you have rotated the knot, it might not
give you the number
you are seeing.  If you do a\\
\hyperlink{rotatefix}{\command{rotate fix}}\\
and then\\
\command{xing}\\
then you will see the number of crossings from the view you are
looking at (at least the \hyperlink{orthographic}{\commandname{orthographic}} version of it).

\hypertarget{parameterssection}{}

\section {Parameter values}
\label{sec:params}

The values in the sliders and rollers are parameters that can be set directly via a command like\\
\hyperlink{sradius}{\command{sradius = 1}}\\
for example. 
This is useful to set a parameter to an exact value, which might be difficult using the slider/roller. 
The majority of parameters in \kp\ do not have a slider/roller associated with them. 
See the \hyperlink{parameters}{\commandname{parameters}} command description.

Some of the most frequently used parameters are the following:

\hypertarget{background}{\commandname{background}} --- Used to change the background colour.  Typically, you do
a \commandname{background = white}.
See the \comm{colour} command for more detail.

\hypertarget{beadrad}{\commandname{bradius}} ---  The radius of the sphere representing the beads (vertices).
Has a roller on \panel{Main}.

\hypertarget{cradius}{\commandname{cradius}} --- This controls the radius of the
cylinders in the knot in  \hyperlink{modecb}{beads and cylinders mode}.  Has a roller on \panel{Main}.

\hypertarget{dbmin}{\commandname{dbmin}} ---  The minimum number of beads that \hyperlink{deletedownto}{\commandname{delete downto}}
will delete to.   Has a slider on \panel{Main}.

\hypertarget{dstep}{\commandname{dstep}} --- The number of relaxation steps between each redisplay. 
Defaults to 1.
Increasing the value helps to speed up relaxations if the rendering time is significant.

\hypertarget{ductrue}{\commandname{duc = true}} --- 
This is good to put at the start of any scripts that are used in production runs.
It will cause \kp\ to have a panic exit if it ``complains'' about anything.
A complaint is made whenever there is a major error, such as a typo causing an unknown command error. 
They are prefaced with the text \verb|***|.
The name comes from ``die upon complaining''.

\hypertarget{hstart}{\commandname{hstart}} , 
\hypertarget{hincr}{\commandname{hincr}} , 
\hypertarget{satur}{\commandname{satur}}  
and \hypertarget{value}{\commandname{value}} --- Controls the hue, saturation, and value (HSV) of knots and links in the default colouring mode (auto).
To get a feel for how these work, go to the \panel{Surf} panel and enter the following commands:\\
\zcomm{reset all}\\
\zcomm{chain 18}\\
and then play with the top four sliders on \panel{Surf}.  

\hypertarget{ncur}{\commandname{ncur}} ---  KnotPlot draws the smooth tubes based on a B\'ezier spline interpolation of
the polygonal knot.  This is the number of interpolation points per edge.
Increase to 10 or more for a smoother look.
Run the \demo{bezier} demo on \panel{DemoB} for an interactive experience.  Has a roller on \panel{Main}.

\hypertarget{nseg}{\commandname{nseg}} ---  The number of subdivisions in the meridional direction for both the smooth tubes and the cylinder.
Load a knot/link and then try the following commands to get a neat extruded triangle look:\\
\zcomm{mode s}\\
\zcomm{nseg = 3}\\
\zcomm{draw hflat}\\
\zcomm{twfix}\\
Has a roller on \panel{Main}.

\hypertarget{N-torus}{\commandname{N-torus}} ---  The default number of beads used if that parameter is omitted for commands 
such as \comm{torus} and \comm{unknot}.

\hypertarget{pserase}{\commandname{pserase}} --- This controls the
amount of white space near crossings in a knot diagram when you use
\hyperlink{psout}{\commandname{psout}}.  See Section \hyperlink{DiagramsSection}{\ref{sec:eps}}.

\hypertarget{psmode}{\commandname{psmode}} --- This controls the output
mode when using \hyperlink{psout}{\commandname{psout}}.  There are
lots of different values.  See Section \hyperlink{DiagramsSection}{\ref{sec:eps}}.

\hypertarget{sradius}{\commandname{sradius}} --- This controls the radius of the
tube about the knot  \hyperlink{modes}{smooth tubes mode}.  Has a roller on \panel{Main}.

\hypertarget{stusplit}{\commandname{stusplit}} --- Sets   how many relaxation steps are done between 
checking for stuck edges.   A value of 0 turns off checking. Has a slider on \panel{Main}.

\hypertarget{vscale}{\commandname{vscale}} --- The scale component of the viewing transform. Has a roller on \panel{Main}.  Note that this affects the view
of the knot, but not the coordinates.

\hypertarget{maintab}{}
\hypertarget{uisection}{}

\section{User interface}
\label{sec:ui}

This section is about buttons, boxes, and rollers in
the \emph{Control Panel} in the \panel{Main} tab.  Below
is a description about what most of the controls on \panel{Main} do.

Everything on this panel can also be done using commands (described
in \zoikt{commands}{Section \ref{sec:commands}}) or by
changing parameter values directly (described in \zoikt{parameterssection}{Section \ref{sec:params}}).

\button{undo} --- undoes last action, equivalent to the command \comm{undo}

\button{reset} --- resets most of KnotPlot's state to the startup values but keeps the currently loaded knot, equivalent to the \comm{reset} command

\button{help} -- a brief description of the options in the selected tab

\button{exit} -- quits \kp{}, equivalent to the commands \comm{quit} and \comm{exit}

\button{All} -- loads a randomly chosen prime knot or link from the Rolfsen catalogue\\
\button{1} -- same as above but restrict to knots (i.e., 1-component links)\\
\button{2} -- restrict to 2-component links\\
\button{3} -- restrict to 3-component links\\
\button{4} -- restrict to 4-component links

\textbf{Knot / Link Zoos}

\button{Rand} -- Randomly chosen from the Rolfsen catalogue\\
\button{A} -- The knots $0_1$ to $9_{28}$\\
\button{B} -- $9_{29}$ through $10_{43}$\\
\button{C} -- $10_{43}$ through $10_{107}$\\
\button{D} -- $10_{108}$ through $10_{165}$ and a few prime 2-component links\\
\button{E} -- loads most of the prime 2-component links\\
\button{F} -- loads prime 2,  3, and 4-component links and one 5-component link\\
\button{Torus} -- loads torus knots and links in order of crossing number; you can keep clicking it to get the next set of more complicated ones\\
\button{RandT} -- loads randomly selected torus knots and links\\

\textbf{Boxes under ``Display'' on the left}

These check boxes toggle various display modes: \pui{smooth} are the
smooth tubes (interpolated from the cylinders/edges), \pui{cylinders}
are centred on the edges of the knot, and \pui{Vcylinders} are like
cylinders but they meet up more nicely.  Edges can also be drawn as
\pui{lines}.  The vertices of the knot can be drawn as \pui{beads}
(spheres centred on the vertices) or \pui{points}.  By selecting
\pui{forces} you will see force vectors showing where the knot is
moving if you press the \pui{go} button.  See the
\zoikt{advanceddynamics}{Section \ref{sec:dynamics}} for more
information about the dynamics.\\

\textbf{In the centre}

\comm{twfix}, \comm{untran}, and \comm{centre} are equivalent to the commands of the same name

\button{smooth tubes} is equivalent to the command \hyperlink{modes}{\commandname{mode s}}

\button{norm} is equivalent to the command \hyperlink{drawnormal}{\commandname{draw normal}}

\button{beads \& cyls} is equivalent to the command \hyperlink{modecb}{\commandname{mode cb}}

\button{spect} is equivalent to the command \hyperlink{drawspectrum}{\commandname{draw spectrum}}

\button{flat} is equivalent to the command \hyperlink{drawflat}{\commandname{draw flat}}

\button{hflat} is equivalent to the command \hyperlink{drawhflat}{\commandname{draw hflat}}\\

\textbf{Boxes on the right}

\button{labels} will let you see the vertex, edge numbers, and component numbers.
The $n^{th}$ edge is the one starting at vertex $n$.  Vertex and edge counting start
at 0.

\button{axes} lets you see the $x$, $y$, and $z$-axes

\button{grid} lets you see the $xy$ plane

\button{orientation} puts some arrows on the knot so that you can see the orientation

\button{backdrop} gives you a coloured background\\

\textbf{Rollers}

\param{vscale} modifies the viewing transform, short for viewing scale

\param{cradius} controls the cylinder radius in \hyperlink{modecb}{\commandname{beads and cylinders}} mode

\hyperlink{beadrad}{\roller{bradius}} controls the bead radius, i.e., the vertex sphere radius, in \hyperlink{modecb}{\commandname{beads and cylinders}} mode

\param{sradius} controls the thickness of the tube in \hyperlink{modes}{\commandname{smooth tubes}} mode

\param{ncur} and \param{nseg} control the level of refinement for the graphics of the cylinders (\param{nseg} only) and 
smooth tubes (both).

\textbf{Dynamics}

Dynamics is a mini version of what you see under the \tabname{Dyna} tab.
For more detailed control of the dynamics, see
\zoikt{advanceddynamics}{Section \ref{sec:dynamics}}.

\button{Init} is the normal dynamics. \button{A}, \button{B}, and
\button{zero} are some other pre-loaded dynamics settings that change
the relaxation.  You can play around with these, but generally you
will not want/need to change these dynamics settings.  The \button{Collide}
box allows the
knot type to change, which  can be useful in some situations.  
Relaxation mode  \button{undamped} should generally be used only when the parameter 
\param{velfo}  (velocity damping, available under the \tabname{Dyna} tab) is turned on.

\roller{dstep} is how many dynamics steps happen between displaying
the new configuration.  If you increase this, things will move faster
because the computer does not need to spend so much time rendering the image.


\param{stusplit} and  \texttt{dbeads}  are usually used in conjunction
to simplify complex or awkward configurations. 
See the exercise in Section \hyperlink{simple}{\ref{sec:simple}} for an example of how to use these features.

\hypertarget{aloadsave}{}

\section{KnotPlot distribution}

\label{sec:aloadsave}

\subsection{Where the files are}

In a similar manner to the  command path inside a terminal shell, KnotPlot has a read, write, and execute path.
Users can see what these are or change the paths using the \comm{path} command. 
On a typical Linux installation, this would print something like
\begin{verbatim}
Current execute path:
   .
   /usr/share/knotplot/resource/linux
   /usr/share/knotplot/resource/generic
Current write path:
   .
Current read path:
   .
   /usr/share/knotplot/basic
   /usr/share/knotplot/special
   /usr/share/knotplot/demos
\end{verbatim}
We can see that the current working directory (represented by \verb|.|) is first in all the paths, so a \comm{load} command
will first look in the that directory for a knot file.  
If it is not found, then the \code{basic}, \code{special}, and \code{demos} directories  are searched in that order.
These directories, as well as \code{resource},  are located under the  KnotPlot home directory\footnote{
\code{/Applications/KnotPlot/KnotPlot.app/Contents/Resources} on macOS and\\
\code{C:\textbackslash Program Files (x86)\textbackslash KnotPlot} on Windows 10/11.\\
In general, use the  command \comm{version} to find out where the KnotPlot home directory is.
}
and have the following function:

\code{basic} --- Contains the Rolfsen catalogue plus a few composite knots.

\begin{figure}
\hbox to \textwidth{\hss\includegraphics[width=4cm]{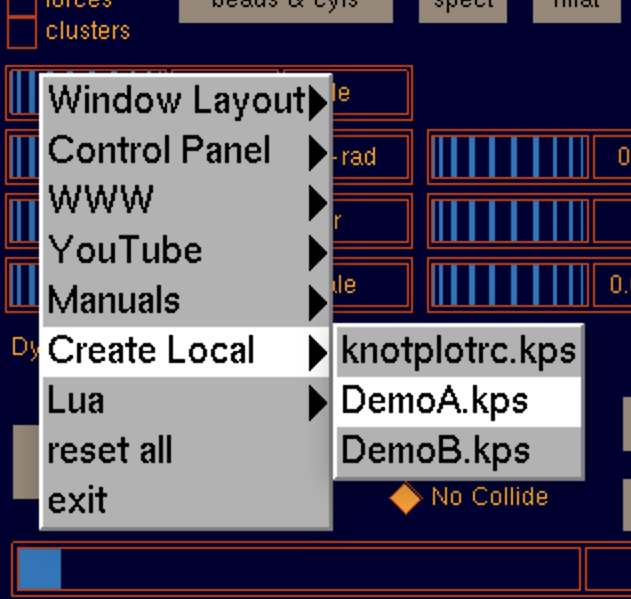}\hss}
  \caption{Creating a starter \code{DemoA.kps} file.}
  \label{fig:local}
\end{figure}

\code{demos} ---  The KnotPlot script files that implement the demos on \panel{DemoA} and \panel{DemoB} are found here. 
Two of these files, \code{DemoA.kps} and \code{DemoB.kps} define which buttons appear on the demo tabs.
Because the working directory is first in the read path, it is possible  to customize the \panel{DemoA} or \panel{DemoB} panels
by creating a local copy of one of these files in that directory.
To create a starter file, right click on the Control Panel and select one of the menu entries as shown in Figure \ref{fig:local}.
The \code{.kps} files are just plain text and can be edited with any text editor. 

\code{resource} ---  A varied collection of resources that KnotPlot relies upon.
The definition for the named colours that the \comm{colour} command recognizes is in the \code{colours.txt} file. 
The \code{src} subdirectory
contains example C++ code to 
read and write KnotPlot's native binary format \cite{KnotPlotFileFormats}  for knots/links.

\code{special} --- A ``Special Collection'' of many knots and links of interest in knot theory and in popular culture.
To load a randomly chosen member of this collection, click on the \pui{spec coll} button on \panel{Main} or enter the command
\code{knot special}.

\subsection{Various catalogues}

As mentioned in Section \hyperlink{LoadNSave}{\ref{sec:loadsave}}
KnotPlot comes with the Rolfsen catalogue in versions that are similar
to Appendix C of Rolfsen's book \cite{Rolfsen} (they are flipped
vertically).  This catalogue is also available in several other
versions, whose knots can be loaded by prepending a subdirectory name
(one of \code{ms}, \code{mseq}, or \code{mscl}) to the knot name.

Minimal stick --- Finding minimal stick representatives for knots is work that was 
pioneered by Richard Randell \cite{RandellEIK}.
RGS extended the known values to all
 prime knot types up to 10 crossings \cite{Scharein_1998}.\\
 \zcomm{load ms/10.124} 

Minimal stick equilateral --- Modifying techniques used in the work  mentioned above, 
the authors \cite{RawSch} found proposed equilateral minimal stick representatives.\\
\zcomm{load mseq/8.19} 

Minimal step cubic lattice   ---  Extending the work by Diao \cite{Diao93}, RGS found minimal step candidates
for all knots and links in the Rolfsen catalogue.  
This work, in addition to some interesting theoretical results, was published in \cite{Scharein_2009}.\\
\zcomm{load mscl/10.123}  

The minimal edge versions of the trefoil are shown in Figure \ref{fig:EPSPL}.

In addition, \kp{} includes the census of the simplest hyperbolic knots by Callahan, Dean, and Weeks \cite{CalDeaWee}.
This study ranks the
 knots arranged according to the complexity of their \emph{complements}, rather than by crossing number as is usually done.
Load the simplest knot in this catalogue with \\
\zcomm{load census/k2.1}\\
and you will see that it is the Figure-8 knot.
The 3D versions of these knots were created by RGS in collaboration with the authors of the study. 

\subsection{Knots and links in Virtual Reality (VR)}

If you happen to be away from KnotPlot, you can still view the knots mentioned in the subsection above as well as knots and links from the Knot Zoo 
by going to one of the web pages\\
\url{https://knotplot.com/zoo/}\\
\url{https://knotplot.com/stick-numbers}\\
\url{https://knotplot.com/hyper/}

On any browser on a computer or mobile device, you will see a \3\ model that you can rotate and scale.
However, if you happen to own a VR headset, you will enter into a space where you can walk around your chosen knot or link and interact with it
in other ways.   

These pages are a work-in-progress and will be enhanced over time. 
Also, example JavaScript code will be supplied in case you want to showcase your own creations in ordinary \3\ or in immersive VR.

\hypertarget{advanceddynamics}{}

\hypertarget{dynatab}{}

\section{Advanced dynamics}     

\label{sec:dynamics}

KnotPlot uses a simple force-based dynamical model to relax and simplify knots.
This model is designed to maintain knot type during relaxation.
An important concept in the model is the notion of \emph{safeness}.
A polygonal knot is said to be in a \emph{safe position} if no two non-adjacent edges are closer than a certain threshold 
distance called \texttt{close}. 
During relaxation vertices are moved one at a time a distance not exceeding \texttt{max-dir}.
Part of the condition for maintaining knot type is that \texttt{max-dir} be strictly less than \texttt{close}.
These parameters can be changed on the \panel{Dyna} tab of the KnotPlot Control Panel or
by setting them directly: \texttt{close = 0.12} and \texttt{max-dir = 0.1}.
You can check to see if your knot is in a safe position using the \comm{safe} command. 
A knot in an unsafe position can usually be made safe by simply scaling it
with the \hyperlink{scale}{\commandname{scale}} command or
to use the \hyperlink{fittomindist}{\commandname{fitto mindist}} command with 
a value larger than \texttt{max-dir}.
For a detailed discussion of KnotPlot's dynamical model and a simple proof that it maintains knot type 
during relaxation see Section 7.1 of RGS's dissertation \cite{Scharein_1998}.

\eric{made the figure bigger}

\begin{figure}[htb] 
\hbox to \textwidth {\hss\includegraphics[height=0.65\textheight]{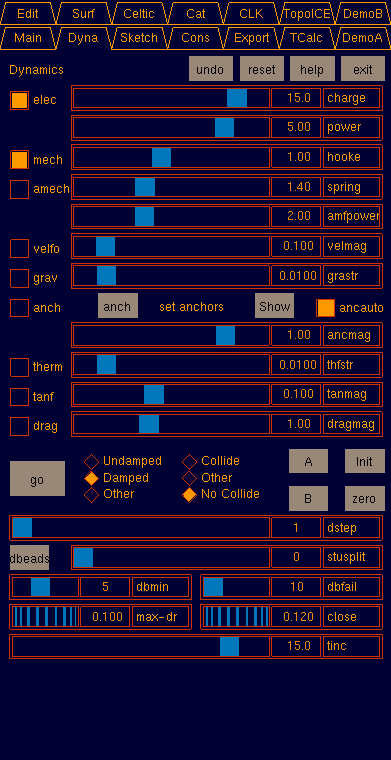}\hss}
\caption{The \panel{Dyna} panel showing options for the dynamics}
\label{fig:Dyna}
\end{figure}

Changes to the dynamical model can be done through the use of commands and parameters 
or by changing settings on the \panel{Dyna} panel.
Figure \ref{fig:Dyna} shows the settings upon startup.
By default only the \texttt{elec}  force (repulsive) and \texttt{mech} force (attractive)  are turned on.  
Forces may be toggled on/off using the check boxes or by typing a command, for example\\
\zcomm{velfo = on}\\
Each of the forces also has a corresponding strength (on the right hand side of the panel)
that can be set using the sliders or by typing a command.\\
\zcomm{velmag = 0.2}\\
At the bottom of this panel are other options, for example checking  the (inappropriately-named) \pui{Collide}  radio button allows 
edges of the knot to pass through each other, thereby allowing the knot type to change. 

As mentioned before, the \pui{dbeads} button deletes as many edges as
it can without changing the knot type.  This is usually done in
conjunction with allowing edges to split if pairs of edges are in a
stuck situation.  Increasing the value of \pui{stusplit} to something
non-zero will cause \kp\ to check every few relaxation steps on
whether there are any stuck edges (according to the value of the
slider).  For a demonstration of the effectiveness of using this
technique in simplifying complex knots, try the \demo{Split/delete}
demo on the \panel{DemoA} control panel.  This demo is essentially
equivalent to the user randomly clicking \pui{dbeads} every so often.

Here is a description of the various forces that are available.
All these forces are applied to the beads/vertices, not the edges.

\pui{elec} ---  A repulsive force similar to an electric field but generally with a higher power in the falloff compared to the
$1/r^2$ 
 Coulomb force.    You can set that \pui{power} as well as the magnitude of the \pui{charge}. 
 This force is applied to all pairs of beads that are not connected by an edge. 
 
 \pui{mech} --- Attractive force similar to a spring but non-linear.
 Change the strength of this force using \pui{hooke}.  The force is
 applied to beads that are connected by an edge.  This allows the user
 to keep some control of the ratio of longest to shortest edge
 lengths.
 
 \pui{amech} --- An alternative to \pui{mech} that attempts to get the edges to a preferred length set by the parameter \pui{spring}. 
 
\hypertarget{velfo}{}

\pui{velfo}  --- A velocity damping force similar to wind resistance.   Only active if the force model is \pui{Undamped}.

\pui{grav}  --- Gravity force.   Try the \demo{knot drop}  and \demo{snakes 2} demos on \panel{DemoB}.

\pui{anch} --- Set some fixed anchor points that are attached to each bead.   The force behaves in a similar way to the \pui{mech} force.
This is good to use if you just want KnotPlot to refine an embedding you have made so that it looks a bit smoother, but not change the configuration too much.
Relaxing a sketched knot using anchors will not generally change the crossing number. 
\eric{ok, I said why}

\pui{therm}  --- A ``thermal'' force that applies a random nudge to each bead at every iteration step.

\pui{tanf} --- Force applied in a tangential direction to each bead in  the direction of its associated edge.
See the \demo{snakes 2} demo where a bunch of  ``snakes'' are dropped onto a plane and they attempt to reptate out of a trap.
Try also the demos \demo{swim} and \demo{dolphins}.  These demos are on \panel{DemoB}.
\eric{added DemoB} 

\pui{drag} --- Enabling this force allows you to click down on a bead and drag it to a new location.  
Try this in in \hyperlink{modecb}{\commandname{beads and cylinders}} mode and run the demo mentioned below.

A good way to get a sense of the various force options is to play with the 
\demo{Forces},  \demo{Anchors}, and \demo{FreeDrag} demos on \panel{DemoA}.

\hypertarget{scripting}{}

\section{Scripting with and running without graphics}
\label{sec:scripting}

Scripting can be used with the \textit{Command Window} or run
in no graphics mode.
For example, suppose you have a bunch of \kp{} commands collected in a plain text file \code{script.kps}.
You can get \kp{} to execute this script by entering \\
\command{<script.kps} \\
into the \textit{Command Window}.

For experiments, \kp{} can be run in non-graphics mode.  
To do this, enter\\
\command{knotplot -nog < script.kps}\\
into a terminal window.
As was mentioned earlier, the command
\comm{ago} (used to run energy dynamics)  is an ``atomic go'', it completes what it is doing before letting anything else happen.
In non-graphics mode, \comm{go} is interpreted  as \comm{ago}. 
In this mode, KnotPlot will not open any windows and commands that require a frame buffer (such as image-producing commands like \comm{imgout}) 
will do nothing. 
However, \comm{psout} will work the same as it  does in graphics mode.

\eric{can't output images}

Through most of KnotPlot's history, scripts have been written in the commands described in this document (and many others).
These \emph{native commands} can do many things, but they are not a general purpose programming language.
In early 2022 such a language, Lua \cite{LuaWikipedia,Lua}, was embedded into \kp.
This addition has  simplified  and greatly enhanced the process of scripting. 

\def\see{\\ \small}
\hypertarget{TabsSection}{}

\section{Tabs}
\label{sec:tabs}

We have only discussed a few of the 14 Control Panel tabs in any detail
(\hyperlink{maintab}{\code{Main}}, 
\hyperlink{xconstab}{\code{Cons}}, 
\hyperlink{sketchtab}{\code{Sketch}}, 
\hyperlink{edittab}{\code{Edit}},
and \hyperlink{dynatab}{\code{Dyna}}).
Here is a brief description of the others.
When you go to any of these tabs, be sure to click on the \pui{help} button at the top of the panel.
The help is specific to the tab you are on.

\pui{Surf} -- Change surface properties of  knots and links (colouring mode, various radii). 
 All of these can  be changed using commands or parameters.
Camera parameters may also be set here such as field of view (FOV) and the distance the camera is from the origin. 
If these sliders are ghosted out, click on the \pui{Perspective} radio button.
Try setting \pui{vcn} to a non-zero value and then adjust \pui{vcA} and \pui{vcphi}.

\pui{Celtic} --- Create Celtic knots.   Quick start: click on one of the 20 preset configurations at the bottom of the panel, for example
\pui{josephine}.   Then click on the \pui{diagram} button at the top.
This will show you the interlace pattern.
Next click on the \pui{copy to arena} button. 
Try the \demo{celtic} demo on \panel{DemoA}.

\pui{Cat} --- Most interesting here is a catalogue of 1736 decorative knots and links with $C_n$ and $D_n$ symmetry for 
$3 \le n \le 8$.

\pui{Export} --- The buttons on this panel are intended as demos only.  
Use the relevant commands for actual work (\comm{psout}, \comm{imgout}, and \comm{export}).
\textbf{Important:} click on the \pui{how to set up a project folder (YouTube)} button on the bottom.

\pui{CLK} --- Knot and links on the simple cubic lattice. 
This panel uses RGS's implementation of the BFACF algorithm \cite{RensWhit91}.

\pui{TCalc} --- An interactive tangle calculator.  
The \comm{tangle} command can be used to send a string of operations to the tangle calculator. 

\pui{TopoICE} --- This is the \emph{Topological Interactive Construction Engine} \cite{TopoICE-X,TopoICE-R}.
It is an interface to software written by Isabel Darcy and uses the tangle calculator extensively. 
This tool is related to site-specific recombination in DNA.
Manuals can be found on the KnotPlot download page \cite{KnotPlotDownload}.

\pui{DemoA}, \pui{DemoB} --- A collection of demos that may be interesting
or even fun.  Briefly, the demos on \panel{DemoA} that most useful for
learning what the dynamics can do are \demo{Split/delete} and
\demo{\pui{Forces}}.  For a complete description, please see the PDF
document (ever evolving) \emph{A Rough Guide to the KnotPlot Demos}
that is located in the KnotPlot installation folder or can be
downloaded from \cite{KPRoughGuide}.

There are several more tabs that can only be accessed by right clicking on the Control Panel and selecting the Control Panel submenu.
Perhaps the most interesting of these is the \panel{4D} tab that allows users to create \4\ knots from \3\ knots.
We will leave that for you to explore.

\hypertarget{ActivitySection}{}

\section{Commands listed by activity}
\label{sec:list}

In this section, We have grouped the commands by certain activities
that you might be trying to do.  The list of activities is below,
then followed by the pertinent commands (whose descriptions are
in Section \ref{sec:commands}).

\subsection{Changing the knot/link coordinates}
\label{changing}
\begin{itemize}
  \item \hyperlink{about}{\commandname{about}}
  \item \hyperlink{nbeads}{\commandname{nbeads}}
  \item \hyperlink{nbeadsmult}{\commandname{nbeads mult}}
  \item \hyperlink{reflect}{\commandname{reflect}}
  \item \hyperlink{revbeads}{\commandname{revbeads}}
  \item \hyperlink{rotatefix}{\commandname{rotate fix}}
  \item \hyperlink{scale}{\commandname{scale}}
  \item \hyperlink{shift}{\commandname{shift}}
  \item \hyperlink{split}{\commandname{split}}
  \item \hyperlink{translate}{\commandname{translate}}
  \item \hyperlink{translateto}{\commandname{translate to}}
\end{itemize}

\subsection{Isolating components in links}
\label{isolating}
\begin{itemize}
  \item \hyperlink{delete}{\commandname{delete}}
  \item \hyperlink{hide}{\commandname{hide}}
  \item \hyperlink{hideall}{\commandname{hide all}}
  \item \hyperlink{keep}{\commandname{keep}}
  \item \hyperlink{unhide}{\commandname{unhide}}
  \item \hyperlink{unhideall}{\commandname{unhide all}}
\end{itemize}

\subsection{Combining multiple knots/links}
\label{combining}
\begin{itemize}
  \item \hyperlink{loadcombine}{\commandname{load combine}}
  \item \hyperlink{loadsum}{\commandname{load sum}}
\end{itemize}

\subsection{Making a knot/link look nice}
\label{nice}
\begin{itemize}
  \item \hyperlink{alignaxes}{\commandname{align axes}}
  \item \hyperlink{centre}{\commandname{centre}}
  \item \hyperlink{centremass}{\commandname{centre mass}}
  \item \hyperlink{deletedownto}{\commandname{delete downto}}
  \item \hyperlink{fitto}{\commandname{fitto}}
  \item \hyperlink{fittomindist}{\commandname{fitto mindist}}
  \item \hyperlink{nbeads}{\commandname{nbeads}}
  \item \hyperlink{nbeadsmult}{\commandname{nbeads mult}}
  \item \hyperlink{rotate}{\commandname{rotate}}
  \item \hyperlink{rotateunit}{\commandname{rotate unit}}
\end{itemize}

\subsection{Making pictures}
\label{pictures}
\begin{itemize}
  \item \hyperlink{background}{\commandname{background}}
  \item \hyperlink{beadrad}{\commandname{bradius}}
  \item \hyperlink{colour}{\commandname{colour}}
  \item \hyperlink{cradius}{\commandname{cradius}}
  \item \hyperlink{dispt}{\commandname{display true}}
  \item \hyperlink{drawflat}{\commandname{draw flat}}
  \item \hyperlink{drawhflat}{\commandname{draw hflat}}
  \item \hyperlink{drawnormal}{\commandname{draw normal}}
  \item \hyperlink{drawspectrum}{\commandname{draw spectrum}}
  \item \hyperlink{imgout}{\commandname{imgout}}
  \item \hyperlink{luxo}{\commandname{luxo}}
  \item \hyperlink{matrgbbead}{\commandname{matrgb bead}}
  \item \hyperlink{matrgbknot}{\commandname{matrgb knot}}
  \item \hyperlink{modecb}{\commandname{mode cb}}
  \item \hyperlink{modes}{\commandname{mode s}}
  \item \hyperlink{orthographic}{\commandname{orthographic}}
  \item \hyperlink{perspective}{\commandname{perspective}}
  \item \hyperlink{psout}{\commandname{psout}}
  \item \hyperlink{rotate}{\commandname{rotate}}
  \item \hyperlink{twfix}{\commandname{twfix}}
  \item \hyperlink{pserase}{\commandname{pserase}}
  \item \hyperlink{psmode}{\commandname{psmode}}
  \item \hyperlink{psoption}{\commandname{psoption}}
\end{itemize}

\subsection{Relaxing a knot/link}
\label{relax}
\begin{itemize}
  \item \hyperlink{ago}{\commandname{ago}}
  \item \hyperlink{deletedownto}{\commandname{delete downto}}
  \item \hyperlink{go}{\commandname{go}}
  \item \hyperlink{massopen}{\commandname{mass open}}
  \item \hyperlink{stusplit}{\commandname{stusplit}}
\end{itemize}

\subsection{Technical issues and information}
\label{technical}
\begin{itemize}
  \item \hyperlink{allocate}{\commandname{allocate}}
  \item \hyperlink{cwd}{\commandname{cwd}}
  \item \hyperlink{ductrue}{\commandname{duc}}
  \item \hyperlink{info}{\commandname{info}}
  \item \hyperlink{jitter}{\commandname{jitter}}
  \item \hyperlink{parameters}{\commandname{parameters}}
  \item \hyperlink{quit}{\commandname{quit}}
  \item \hyperlink{reset}{\commandname{reset}}
  \item \hyperlink{resetall}{\commandname{resetall}}
  \item \hyperlink{undo}{\commandname{undo}}
  \item \hyperlink{untran}{\commandname{untran}}
\end{itemize}

\subsection{Knot/Link configuration measurements}
\label{measurements}
\begin{itemize}
  \item \hyperlink{acn}{\commandname{acn}}
  \item \hyperlink{angle}{\commandname{angle}}
  \item \hyperlink{energy}{\commandname{energy}}
  \item \hyperlink{length}{\commandname{length}}
  \item \hyperlink{lnknum}{\commandname{lnknum}}
  \item \hyperlink{rog}{\commandname{rog}}
  \item \hyperlink{thickness}{\commandname{thickness}}
\end{itemize}

\subsection{Scripting}
\label{scripting}
\begin{itemize}
  \item \hyperlink{ago}{\commandname{ago}}
  \item \hyperlink{lua}{\commandname{lua}}
  \item \hyperlink{modecb}{\commandname{mode cb}}
  \item \hyperlink{modes}{\commandname{mode s}}
  \item \hyperlink{nap}{\commandname{nap}}
\end{itemize}

\subsection{Crossing information}
\label{crossinginfo}
\begin{itemize}
  \item \hyperlink{diagram}{\commandname{diagram}}
  \item \hyperlink{dowker}{\commandname{dowker}}
  \item \hyperlink{gauss}{\commandname{gauss}}
  \item \hyperlink{xing}{\commandname{xing}}
\end{itemize}

\subsection{Special knot constructions}
\label{constructs}
\begin{itemize}
\item \comm{braid}
  \item \hyperlink{torus}{\commandname{torus}}
  \item \hyperlink{unknot}{\commandname{unknot}}
\end{itemize}

\subsection{Safety}
\label{safety}
\begin{itemize}
  \item \hyperlink{distance}{\commandname{distance}}
  \item \hyperlink{fittomindist}{\commandname{fitto mindist}}
  \item \hyperlink{safe}{\commandname{safe}}
\end{itemize}

\subsection{File input/output}
\label{fileio}
\begin{itemize}
  \item \hyperlink{load}{\commandname{load}}
  \item \hyperlink{loadcombine}{\commandname{load combine}}
  \item \hyperlink{loadsum}{\commandname{load sum}}
  \item \hyperlink{save}{\commandname{save}}
\end{itemize}

\bibliographystyle{alpha}
\bibliography{knotplot}
\end{document}